\documentclass[a4paper,twoside,
11pt]{amsart}
\usepackage{amssymb,amsfonts,amsmath}
\begin{document}
\newtheorem{defin}{~~~~Definition}%[section]
\newtheorem{prop}{~~~~Proposition}[section]
\newtheorem{remark}{~~~~Remark}[section]
\newtheorem{cor}{~~~~Corollary}%[section]
\newtheorem{theor}{~~~~Theorem}%[section]
\newtheorem{lemma}{~~~~Lemma}[section]
\newtheorem{ass}{~~~~Assumption}
\newtheorem{con}{~~~~Conjecture}
\newtheorem{concl}{~~~~Conclusion}
\renewcommand{\theequation}{\thesection.\arabic{equation}}
\newcommand {\trans} {^{\,\mid\!\!\!\cap}}
\newcommand{\vf}{\varphi}
\newcommand{\mJ}{\mathcal J}
\newcommand{\e}{\varepsilon}
\title{ On local geometry of nonholonomic rank 2 distributions
% A Canonical Frame for Nonholonomic
%Rank Two Distributions of Maximal Class
}
\author
{Boris Doubrov
%$^1$ \footnote
\address{
%$~^1$
Belarussian State University, Nezavisimosti Ave.~4, Minsk 220050,
Belarus;
%\email {
 E-mail: doubrov@islc.org}\and Igor Zelenko
%$^2$
\address
%\footnote
{
%$~^2$
S.I.S.S.A., Via Beirut 2-4, 34014, Trieste,
Italy;
%}
%}
%\email {
E-mail:
 zelenko@sissa.it}}
 \maketitle\markboth{Boris Doubrov and Igor Zelenko} {On
Local Geometry of Rank 2 Distributions}

\begin{abstract}
In 1910 E. Cartan constructed a canonical frame and found the most
symmetric case for maximally nonholonomic rank $2$ distributions in
$\mathbb R^5$.
%, applying his method of equivalence.
We
%introduce the notion of the class of completely nonalcoholic rank
%$2$ distributions in ${\mathbb R}^n$ and
solve the analogous problem for germs of generic rank $2$
distributions in ${\mathbb R}^n$ for $n>5$. We use a completely
different approach based on the symplectification of the problem.
The main idea is to consider a special odd-dimensional submanifold
$W_D$ of the cotangent bundle associated with any rank 2
distribution $D$. It is naturally foliated by characteristic curves,
which are also called the abnormal extremals of the distribution
$D$. The dynamics of vertical fibers along characteristic curves
defines certain curves of flags of isotropic and coisotropic
subspaces in a linear symplectic space. Using the classical theory
of curves in projective spaces, we construct the canonical frame of
the distribution $D$ on a certain $(2n-1)$-dimensional fiber bundle
over $W_D$ with the structure group of all M\"obius transformations,
preserving $0$. The paper is the detailed exposition of the
constructions and the results, announced in the short
note~\cite{doubzel}.
\end{abstract}

\section{Introduction}
\setcounter{equation}{0}
 A
rank $l$ vector distribution $D$ on an $n$-dimensional
manifold $M$ or an $(l,n)$-distribution (where $l<n$) is a
subbundle of the tangent bundle $TM$ with $l$-dimensional
fibers. The group of germs of diffeomorphisms of $M$ acts
naturally on the set of germs of $(l,n)$-distributions and
defines the equivalence relation there.
%on this set.
%In other words, for each point $q\in M$ a
%$l$-dimensional subspace $D(q)$ of the tangent space $T_qM$
%is chosen and $D(q)$ depends smoothly on $q$.
%The group of germs of diffeomorphisms at a point $q\in M$
%acts naturally on the germs of rank $l$ distributions.
%Two vector distributions $D_1$ and $D_2$ are called
%equivalent, if there exists a diffeomorphism $F:M\mapsto M$
%such that
%$F_*D_1(q)=D_2(F(q))$ for any $q\in M$.
% Two germs of vector
%distributions $D_1$ and $D_2$ at the point $q_0\in M$ are
%called equivalent, if there exist neighborhoods $U$ and
%$\tilde U$ of $q_0$ and a diffeomorphism $F:U\mapsto \tilde
%U$ such that $F_*D_1(q)=D_2(F(q))$ for all $q\in U$ and
%$F(q_0)=q_0$.
The question is \emph{when two germs of distributions are
equivalent?} Distributions are naturally associated with
Pfaffian systems and with control systems linear in the
control. So the problem of equivalence of distributions can
be reformulated as the problem of equivalence of the
corresponding Pfaffian systems and the state-feedback
equivalence of the corresponding control systems. The
obvious (but very rough in the most cases) discrete
invariant of a distribution $D$ at $q$ is so-called \emph{
the small growth vectors} at $q$. It is the tuple
%$\bigl(\dim D(q),\dim D^2(q),\dim D^3(q),\ldots\bigr)$,
$\{\dim D^j(q)\}_{j\in{\mathbb N}}$, where $D^j$ is the
$j$-th power of the distribution $D$, i.e.,
$D^j=D^{j-1}+[D,D^{j-1}]$, $D^1=D$. %, %$j\geq 2$.
% the big
%growth vector is
%the tuple $\bigl(\dim D(q),\dim
%D_2(q),\\\dim D_3(q),\ldots\bigr)$, where $D_j$ is defined by the
%following recursive formula
%$D_j=D_{j-1}+[D_{j-1},D_{j-1}]$, $D_1=D$.
%The basic question is whether two given germs of
%distributions lie in the same orbit w.r.t. this action.
%Such germs of distributions are called equivalent.
A simple estimation shows that at least $l(n-l)-n$
functions of $n$ variables are required to describe generic
germs of $(l,n)$-distribution, up to the equivalence (see
\cite{versh1} and \cite {zhit0} for precise statements).
 There are
only three cases, where $l(n-l)-n$ is not positive: $l=1$
(line distributions), $l=n-1$, and $(l,n)=(2,4)$. Moreover,
it is well known that in these cases generic germs of
distributions are equivalent.
%It is well known that generic germs of distributions are
%equivalent only in the following three cases: $l=1$,
%%(line distributions),
%$l=n-1$, and $(l,n)=(2,4)$.
%:
%%Darboux's
%%theorem in the first case and
%%Engel's theorem in the second case,
For $l=1$ it is just the classical theorem about the
rectification of vector fields without stationary points,
for $l=n-1$ all generic germs are equivalent to Darboux's
model, while for $(l,n)=(2,4)$ they are equivalent to
Engel's model (see, for example, \cite{bryantbook}). In all
other cases generic $(l,n)$-distributions have functional
invariants.

In the present paper we restrict ourselves to the case of rank 2
distributions, although our method can be applied also for
distribution of rank greater than 2 (as will be described in the
forthcoming paper \cite{doubzel3}). The model examples of rank 2
distributions come from so-called underdetermined ODE's of the type
\begin{equation}
\label{under}
%\[
z'(x)=F\bigl(x,y(x),\ldots, y^{(n-3)}(x),z(x)\bigr),
%\quad
%r=n-3,
%\]
\end{equation}
for two functions $y(x)$ and $z(x)$. Setting $p_i=y^{(i)}$,
$0\leq i\leq n-3$, with each such equation one can associate the rank $2$
distribution
%${\mathcal D}$
%system of
%$n-2$ Pfaffian equations
on ${\mathbb R}^n$ with coordinates $(x,p_0,\ldots,
p_r,z)$ defined as the annihilator of the following $n-2$ 1-forms:
%\begin{equation*}
%\label{Pfaff}
 %\left\{\begin{array}{l}
\begin{equation}
\label{Pfaff}
\begin{aligned}
~&dp_i-p_{i+1} dx , \,\,0\leq i\leq n-4,\\
 ~&dz-F(x,p_0,\ldots, p_{n-3},z)dx.
\end{aligned}
\end{equation}
%\end{array}
%\right.,
%\end{equation*}
%which in turn defines rank 2 distribution in ${\mathbb
%R}^n$.

%{\bf 2. The class of a rank 2 distribution.}
%
%{\bf 3. Jacobi curves and the main theorem}
%
%{\bf 4. Proof of Theorem \ref{main}}
%
%{\bf 5. Appendix A. The projective structures on the curves
%in Grassmannians of half-dimensional subspaces}
%
%{\bf 6. Appendix B. On application of the Tanaka-Morimoto
%theory}

%In the present paper we are
%interested in finding of functional invariants of
%distributions w.r.t. the action of the group of
%diffeomorphisms.
% w.r.t. the action of the group of the
%diffeomorphisms, which help in solving the equivalence
%problem.

%The goal of the present paper is to construct basic
%functional (differential) invariants of distributions. Using these invariants o%ne  can say if two germs of distributions are
%equivalent or not.

%Let us roughly estimate the "number of parameters" in the
%considered equivalence problem. The set of $l$-dimensional
%subspaces in $\mathbb{R}^n$ forms $l(n-l)$-dimensional
%manifold. Therefore, if the coordinates on $M$ are fixed
%then the rank $l$ distribution can be defined by $l(n-l)$
%functions of $n$ variables. The group of the coordinate
%changes
%on $M$ is parameterized by $n$ functions of $n$ variables.
%So, by a coordinate change one can "normalize", in general,
%only $n$ functions among those $l(n-l)$ functions, defining
%the distribution.
%It is well known that the set of
%classes of equivalent germs of rank $l$ distributions
%can be parameterized by $l(n-l)-n$ arbitrary germs of
%functions of $n$ variables .

For $n=3$ and $4$ all generic germs of rank 2 distribution
are equivalent to the distribution, associated with the
underdetermined ODE $z'(x)=y(x)$ (Darboux and Engel models
respectively).
%(the reformulaion
%of Darboux's and Engel's theorems respectively).
The case
%$l=2$,
$n=5$ (the smallest dimension, when functional parameters
appear)
%functional invariant should
%appear starting from $n=5$.
%It is well
%known that in the low dimensions $n=3$ or $4$ all generic
%germs of rank 2 distributions are equivalent. .
%The case of generic $(2,5)$-distributions
was treated by E. Cartan in \cite{cartan} with his
reduction-prolongation procedure. First, for any
$(2,5)$-distribution with the small growth vector $(2,3,5)$
he constructed the canonical coframe in some 14-dimensional
manifold, which implied that the group of symmetries of
such distributions is at most $14$-dimensional. Second,
%E. Cartan
he showed that any $(2,5)$-distribution with
$14$-dimensional group of symmetries is locally equivalent
to the distribution, associated with the underdetermined
ODE $z'(x)=\bigl(y''(x)\bigr)^2$, and its group of
symmetries
%of this distribution
is isomorphic to the real
split form of the exceptional Lie group $G_2$. Historically
it was the first natural appearance of this group.

After the work of Cartan the open question was \textbf {to
construct the canonical frame and to find the most
symmetric cases for $(2,n)$-distributions with $n>5$}.
%Actually, Cartan's application of his equivalence method to
%$(2,5)$-distributions is the most complicated application
%of this method till now.
% When the dimension of the ambient
%manifold increases, the computational difficulties in
%applying of Cartan's method increase as well.
The Cartan equivalence method was systematized and
generalized by N. Tanaka and T. Morimoto (see \cite{tan2,
mori}). Their theory is heavily based on the notion of
so-called \emph{ symbol algebra} of the distribution at a
point, which is a special \emph{ graded nilpotent Lie
algebra}, naturally associated with the distribution at a
point. The symbol algebras have to be isomorphic at
different points and all constructions strongly depend on
the type of the symbol. Note that already in the case of
$(2,6)$-distributions with maximal possible small growth
vector $(2,3,5,6)$ three different symbol algebras are
possible, while for $n=9$ the set of all possible symbol
algebras depends on continuous parameters, which implies in
particular that generic distributions do not have a
constant symbol.
%(see, for example, \cite{marigo}).

In the present paper we give the answer to the question, underlined
in the previous paragraph, for rank $2$ distributions from some
generic class. Our constructions are based on a completely
different, variational approach, developed in \cite{jac1} and
\cite{zelvar}. The paper is  the detailed exposition of the
constructions and the results, announced in short note
\cite{doubzel}.
%This approach was applied
%to the construction of differential invariants of rank 2
%distributions in \cite{zelvar}. In the present paper we use
%this approach to
%we give the answer to the question, underlined in the
%previous paragraph, for
%rank $2$
%distributions from some generic class.
%Roughly speaking, we make a kind of symplectification of
%the problem by lifting the distribution to the cotangent
%bundle $T^*M$ of the
%ambient
%manifold $M$.
The starting point of this approach is to lift a distribution $D$ to
a special odd-dimensional submanifold $W_D$ of the cotangent bundle,
foliated by the characteristic curves, which are also called the
abnormal extremals of the distribution $D$. They are Pontryagin
extremals with zero Lagrange multiplier near the functional for any
extremal problem with constrains, given by the distribution $D$.
%while if $D$ is a distribution of odd rank it is the annihilator of
%$D$ itself.
%, denoted by $D^\perp$.
The dynamics of the lifting (to $W_D$)
%$(D^2)^\perp$)
of the distribution $D$ w.r.t. to this $1$-foliation along
any abnormal extremal
%of $D$
%characteristic curve (of  $W_D$)
%$(D^2)^\perp$)
defines certain
%unparameterized
curve of flags of isotropic and coisotropic subspaces in a linear
symplectic space. So, the problem of equivalence of distributions
can be essentially reduced to the differential geometry of such
curves: symplectic invariants of these curves automatically produce
invariants of the distribution $D$ itself and the canonical frame bundles
associated with such curves can be in many cases effectively used
to construct the canonical frames of the distribution
$D$ itself defined on a certain bundle over $W_D$.

In the case of a nonholonomic rank 2 distribution  the submanifold
$W_D$ is nothing but the annihilator of the square of $D$, denoted
by $(D^2)^\perp$. Under additional generic assumptions the curves of
flags, associated with abnormal extremals from a certain open and
dense subset of  $(D^2)^\perp$, are curves of complete flags.
Moreover, these curves of complete flags can be recovered by
differentiation from the curves of their one-dimensional subspaces,
i.e. from  curves in projective spaces. Recall that the differential
geometry of curves in projective spaces was developed already in
1905 by E.J.~Wilczynski (\cite{wil}). In particular, these curves
(and therefore the corresponding abnormal extremals of the
distribution) are endowed with the canonical projective structure,
i.e., we have a distinguished set of parameterizations (called
projective) such that the transition function from one such
parameterization to another is a M\"{o}bius transformation. Besides,
for each fixed projective parameterization on such curve one can
construct the canonical moving frame in the ambient linear
symplectic space.
%These curves are endowed with the canonical projective structure and
%given a parametrization on them one can find the moving frame,
%canonically associated with this parametrization.

These two facts together allow us to construct the canonical frame
for any $(2,n)$-distribution $D$, $n>5$, from a certain generic
class. This frame lives on a certain principle bundle over
$(D^2)^\perp$ with the structure group $ST(2,\mathbb R)$ of all
M\"obius transformations, preserving $0$. The fiber of this bundle
over the point $\lambda\in (D^2)^\perp$ is just the set of all
projective parameterizations of the abnormal extremal passing
through $\lambda$ such that the point $\lambda$ corresponds to $0$.
In particular, it implies that the group of symmetries of such
distributions is at most $(2n-1)$-dimensional.

\section {Abnormal extremals}
%The class of rank 2 distribution}
\setcounter{equation}{0}

Assume that $\dim D^2(q)=3$ and $\dim D^3(q)>3$ for any
$q\in M$. Denote by $(D^l)^{\perp}\subset T^*M$ the
annihilator of the $l$th power $D^l$, namely
%\begin{equation}
%\label{annihil}
$$(D^l)^{\perp}= \{(p,q)\in T^*M:\,\, p\cdot v=0\,\,\forall
v\in D^l(q)\}.$$
%\end{equation}
First we distinguish a characteristic $1$-foliation on the
codimension $3$ submanifold
$(D^2)^\perp\backslash(D^3)^\perp$ of $T^*M$. For this let
$\pi:T^*M\mapsto M$ be the canonical projection. For any
$\lambda\in T^*M$, $\lambda=(p,q)$, $q\in M$, $p\in
T_q^*M$, let $\mathfrak{s}(\lambda)(\cdot)=p(\pi_*\cdot)$
be the canonical Liouville form and $\sigma=d\mathfrak {s}$
be the standard symplectic structure on $T^*M$. Since the
submanifold $(D^2)^\perp$ has odd codimension in $T^*M$,
%Therefore
the kernels of the restriction $\sigma|_{(D^2)^\perp}$ of $\sigma$
on $(D^2)^\perp$ are not trivial. Moreover, as we show below, for
%the generic points of
%$(D^2)^\perp$, namely,
the points of $(D^2)^\perp\backslash (D^3)^\perp$ these
kernels are one-dimensional.
%These kernels
They form the \emph{characteristic line distribution} in $(D^2)^\perp\backslash(D^3)^\perp$, which will be denoted by
${\mathcal C}$. The line distribution ${\mathcal C}$ defines a \emph{characteristic 1-foliation} of
$(D^2)^\perp\backslash(D^3)^\perp$.
%The leaves of this
%foliation are called  the {\it characteristic curves} of
%$(D^2)^\perp\backslash(D^3)^\perp$.
 %In Control Theory
Actually the leaves of this
foliation
%these characteristic
%curves
are
%also
so-called \emph{regular abnormal extremals of the distribution $D$}.

Recall that abnormal extremals of $D$ are by definition Pontryagin extremals with zero Lagrange multiplier near the
functional for any extremal problem with constrains, given by the distribution $D$ (and  so they depend only on $D$ and
not on a functional). Regularity means that they do not pass through $(D^3)^\perp$, which is equivalent to the fact
that they satisfy so-called strong generalized Legendre--Glebsch condition (\cite{agrsach}, \cite{zel}). In the sequel
for shortness we will omit the word regular and the leaves of the characteristic foliation will be called just \emph
{abnormal extremals of $D$}.

Let us describe the characteristic line distribution
${\mathcal C}$ in terms of a local basis $(X_1, X_2)$ of
the distribution $D$,
%\begin{equation}
%\label{X12}
$D( q)={\rm span}\{X_1(q), X_2(q)\}$.
%\end{equation}
%Since our study is local, we can always suppose that such
%basis exists, restricting ourselves, if necessary, on some
%coordinate neighborhood instead of whole $M$.
%Given the
%basis $X_1$, $X_2$ one can construct a special vector field
%tangent to the characteristic $1$-foliation $Ab_D$. For
%this suppose that
%\begin{eqnarray}
Denote by
\begin{equation}
\label{x345}
%&~&
X_3=[X_1,X_2],\,\,
%\quad {\rm mod}\, D,\,\,\,
X_4 =\bigl[X_1,[X_1,X_2]\bigr],\,\,
%=[X_1,X_3]\quad {\rm mod}\, D^2,
%\nonumber
%\\ &~&~
%\label{x345} \\
%&~&
X_5 =\bigl[X_2,[X_1,X_2]\bigr].
%=[X_2,X_3]\quad{\rm mod}\, D^2.
%\nonumber
%\end{eqnarray}
\end{equation}
Let us introduce the ``quasi-impulses''
$u_i:T^*M\mapsto\mathbb R$, $1\leq i\leq 5$,
\begin{equation}
\label{quasi25} u_i(\lambda)=p\cdot
X_i(q),\,\,\lambda=(p,q),\,\, q\in M,\,\, p\in T_q^* M
\end{equation}
Then by definitions
\begin{equation}
\label{d2u}
(D^2)^\perp=\{\lambda\in T^*M:
u_1(\lambda)=u_2(\lambda)=u_3(\lambda)=0\}.
\end{equation}
 As usual, for
given function $G:T^*M\mapsto \mathbb R$ denote by $\vec G$
the corresponding Hamiltonian vector field defined by the
relation $i_{\vec G}\sigma
%(\vec G,\cdot)
=-d\,G
%(\cdot)
$.

\begin{lemma}
%Then it is
%easy to show (see, for example, \cite{zel}) that
%\begin{equation}
%\label{ker25} \ker\sigma\Bigr|_{D^\perp}\Bigl.(\lambda)=
%{\rm span}(\vec u_1(\lambda),\vec u_2(\lambda)),\quad
%\forall \lambda\in D^\perp,
%\end{equation}
The characteristic line distribution $\mathcal C$ satisfies
\begin{equation}
\label{foli25}
%\ker\sigma\Bigr|_{(D^2)^\perp}\Bigl.(\lambda)
{\mathcal C}= \langle u_4\vec{u}_2-u_5\vec{u}_1 \rangle.
%{\mathcal C}= \bigl\{\mathbb{R} \bigl((u_4
%\vec{u}_2-u_5\vec{u}_1)\bigr)\bigr\}.
%\quad \forall
%\lambda\in (D^2)^\perp\backslash (D^3)^\perp
\end{equation}
\end{lemma}

\begin{proof}
Take a vector field $H$ on $(D^2)^\perp\backslash(D^3)^\perp$ such
that locally ${\mathcal C}(\lambda)= \bigl\{\mathbb{R}
H(\lambda)\}$. Then by definition of $\mathcal C$ we have
$i_H\sigma|_{(D^2)^\perp}=0$. From this and \eqref{d2u} it follows
that $i_H\sigma \in
%{\rm span} \{d\,u_i\}_{i=1}^3$
\langle du_1, du_2, du_3 \rangle,$ which implies that
\begin{equation}
\label{Hsp} H\in \langle \vec u_1, \vec u_2, \vec u_3 \rangle.
%H\in {\rm span} \{\vec u_i\}_{i=1}^3.
\end{equation}
On the other hand, $H$ is tangent to $(D^2)^\perp$, i.e $du_j(H)=0$ for $1\leq j\leq 3$. This and  (\ref{Hsp}) easily
implies (\ref{foli25}).
\end{proof}
% (or
%{\it the abnormal extremals}) of the distribution $D$ (the
%second term comes from Optimal Control Theory).
% these
%characteristic curves are also called {\it regular abnormal
%extremals of $D$}.
%As a consequence, if $\dim D^3(q)=5$, then
As a consequence, we see that $\pi_*({\mathcal C}(\lambda))\subset
D(q)$ for any $\lambda \in (D^2)^\perp\backslash (D^3)^\perp$,
$\pi(\lambda)=q$. Moreover, if $\dim D^3(q)=5$, then one-dimensional
subspaces $\pi_*({\mathcal C}(\lambda))$ generate $D(q)$:
\begin{equation}
\label{projC} {\rm span}\left\{\pi_*({\mathcal C}(\lambda))\,:\,\lambda\in (D^2)^\perp\backslash
(D^3)^\perp,\pi(\lambda) =q\right\}
%\Bigr)\Bigr)
=D(q).
\end{equation}
In particular, in this case the original distribution can
be recovered from its characteristic line distribution.

% In the sequel given two
%submanifold $S_1$ and $S_2$ of the tangent bundle of some
%manifold $W$ such that $S_i(w)=S_i\cap T_w W$, $i=1,2$, are
%linear subspaces of $T_wW$ (not necessary of the same
%dimensions for different $w$) we will denote by $[S_1,
%S_2]$ the subset $\{[S_1, S_2](w)\}_{w\in W}$ of $TW$ such
%that $$[S_1, S_2](w)={\rm span}\{[Z_1, Z_2](w):Z_i\,\,{\rm
%are}\,\, {\rm vector}\,\, {\rm fields}\,\, {\rm
%tangent}\,\,{\rm to}\,\, S_i, i=1,2\}.$$ It is easy to show
%that with such definition $S_i\subset [S_1, S_2]$ ,
%$i=1,2$.

%Secondly to any point $q\in M$ we assign a natural number,
%%an integer number between $1$ and $n-3$,
%called {\it the class
%of D at $q$}. For this let

\section{The curves of flags associated with abnormal extremals}
\setcounter{equation}{0}

 Now, following \cite{zelvar}, let $\mJ$ be
the pull-back of the distribution $D$ on
$(D^2)^\perp\backslash(D^3)^\perp$ by the canonical projection
$\pi$:
\begin{equation}
\label{prejac}
{\mathcal J}(\lambda)=
%\bigl(T_\lambda
%(T^*_{\pi(\lambda)}M)+
%\ker\sigma|_{D^\perp}(\lambda)\bigr)\cap T_\lambda
%(D^2)^\perp=
%\begin{equation}
%\label{jacproj}
%\hat L_\T=
\{v\in T_{\lambda}(D^2)^\perp:\,\pi_*\,v\in D(\pi\bigl(\lambda)\bigr)\}.
\end{equation}
%\end{equation}
%(\begin{scriptsize}here $T_\lambda (T^*_{\pi(\lambda)}M)$ is the tangent space to the
%fiber $T^*_{\pi(\lambda)}M$ at the point $\lambda$).\begin{footnotesize}\begin{small}\begin{normalsize}\begin{large}\begin{Large}\begin{LARGE}\begin{huge}\end{huge}\end{LARGE}\end{Large}\end{large}\end{normalsize}\end{small}\end{footnotesize}\end{scriptsize}
%Note that $\dim {\mathcal J}(\lambda)=n-1$.
%Actually, ${\mathcal J}$ is
Note that $\dim \mJ = n-1$ and $\mathcal C\subset \mJ$ by \eqref{foli25} .
%The rank $n-1$
%distribution ${\mathcal J}$ is called the \emph {lift of
%distribution $D$ to $(D^2)^\perp\backslash(D^3)^\perp$}. We will
%work with the lift $\mJ$ instead of the original distribution $D$.
The distribution ${\mathcal J}$ is called the \emph {lift of
distribution $D$ to $(D^2)^\perp\backslash(D^3)^\perp$}.

In the sequel we shall work with the lift $\mJ$ instead of the original distribution~$D$. The crucial advantage of
working with $\mJ$ is that it has the distinguished line sub-distribution $\mathcal C$, while the original distribution
$D$ has no distinguished sub-distributions in general.

We can produce two monotonic (by inclusion) sequences of
distributions (in general of nonconstant ranks) first by making
iterative Lie brackets of $\mathcal C$ and $\mJ$ and then by taking
skew symmetric complements w.r.t. the form $\sigma$ of the subspaces
obtained in the previous step. Namely, first define a sequence of
subspaces ${\mathcal J}^{(i)}(\lambda)$, $\lambda\in
(D^2)^\perp\backslash (D^3)^\perp$, by the following recursive
formulas:
\begin{equation}
\label{Ji1} {\mathcal J}^{(i)}= {\mathcal J}^{(i-1)}+ [{\mathcal C},{\mathcal J}^{(i-1)}], \quad {\mathcal
J}^{(0)}={\mathcal J},
\end{equation}
and then set
\begin{equation}
\label{contr} {\mathcal J}_{(i)}(\lambda)= \{v\in T_\lambda
\bigl((D^2)^\perp\bigr): \sigma (v, w)=0\,\, \forall w\in {\mathcal
J}^{(i)}(\lambda)\},
\end{equation}

 We summarize the main properties of  the sequences
$\{{\mathcal J}^{(i)}\}_{i\geq 0}$  and $\{{\mathcal
J}_{(i)}\}_{i\geq 0}$ in the following:
\medskip
\begin{prop}
\label{proper} ~

\begin{enumerate}
\item $\sigma|_\mJ=0$,  $\mJ_{(0)}=\mJ^{(0)}$;

\item $\mJ^{(i-1)}(\lambda)\subseteq {\mathcal J}^{(i)}(\lambda)$,
 $\mJ_{(i)}(\lambda)\subseteq {\mathcal J}_{(i-1)}(\lambda)$;

\item ${\rm dim}\,{\mathcal J}^{(1)}(\lambda)-{\rm dim}\, {\mathcal J}(\lambda) =1$,  ${\rm dim}\,{\mathcal
J}(\lambda)-{\rm dim}\, {\mathcal J}_{(1)}(\lambda) =1$;

\item ${\rm dim}\,{\mathcal J}^{(i)}(\lambda)-{\rm dim}\, {\mathcal J}^{(i-1)}(\lambda)\leq 1$, ${\rm dim}\,{\mathcal
J}_{(i-1)}(\lambda)-{\rm dim}\, {\mathcal J}_{(i)}(\lambda)\leq 1$ for $i\ge2$;

\item ${\rm dim}\, {\mathcal J}^{(i)}(\lambda)\leq 2n-4$.
%\item$\{\mathcal C\}, e\}
\end{enumerate}
\end{prop}

\begin{proof} First note that the second relations in the Properties
(1)-(4) are direct consequences of the corresponding first
relations.

Now let us prove the first relation in the Property (1). By the
arguments similar to the proof of Lemma \ref{foli25}, the set of
points, where the form $\sigma|_{D^\perp}$ is degenerated, coincides
with $(D^2)^\perp$ and for each $\lambda\in (D^2)^\perp$ the kernel
of $\sigma|_{D^\perp}(\lambda)$ satisfies
\begin{equation}
\label{kerDperp} {\rm ker}\, \sigma|_{D^\perp}(\lambda)={\rm
span}\{\vec u_1(\lambda),\vec u_2(\lambda)\},
\end{equation}
where $u_i$ are as in  \eqref{quasi25} for some local basis
$(X_1,X_2)$ of the distribution $D$. Also denote by $\widetilde
V(\lambda)$ the vertical subspace of $T_\lambda D^\perp$, i.e.
$\widetilde V(\lambda)=\{v\in T_\lambda D^\perp,\pi_*v=0\}$.  Then
from \eqref{prejac} and \eqref{kerDperp} one gets easily that
$$\mJ(\lambda)=\bigl(\widetilde V(\lambda)+{\rm ker}\, \sigma|_{D^\perp}(\lambda)\bigr)\cap T_{\lambda}(D^2)^\perp.$$
This immediately implies the first relation of Property~(1).
%from which the first relation of Property (1) follows immediately.
%here $T_\lambda (T^*_{\pi(\lambda)}M)$ is the tangent space to the
%fiber $T^*_{\pi(\lambda)}M$ at the point $\lambda$. Since for any $\lambda\in (D^2)^\perp$  and $v\in

The first inclusion in Property (2) follows directly from definition
\eqref{Ji1} of ${\mathcal J}^{(i)}$. Further, one can easily get
that
\begin{equation} \label{J1} {\mathcal
J}^{(1)}(\lambda)=\{v\in T_{\lambda}(D^2)^\perp:\,\pi_*\,v\in
D^2(\pi\bigl(\lambda)\bigr)\}.
\end{equation}
Since $\dim D^2-\dim D=1$, we obtain the first relation in Property
(3). Besides, since on each step of the recursive formula
\eqref{Ji1} one makes Lie brackets with the rank 1 distribution
$\mathcal C$, Property (3) immediately implies Property (4). In
order to prove Property (5) note that the line distribution
${\mathcal C}$ forms the Cauchy characteristic of the corank $1$
distribution on $(D^2)^\perp$, given by the Pfaffian equation
$\mathfrak {s}|_{(D^2)^\perp}=0$, where as before $\mathfrak {s}$ is
the Liouville form. Since by construction $${\mathcal
J}\subset\{{\mathfrak s}|_{(D^2)^\perp}=0\},$$ one has
\begin{equation}
\label{inLiuv} {\mathcal J}^{(i)}\subset\{{\mathfrak
s}|_{(D^2)^\perp}=0\}\quad \text{for all } i\in\mathbb N.
\end{equation}
Property (5)  follows from the fact that the distribution
$\{{\mathfrak s}|_{(D^2)^\perp}=0\}$ has rank $2n-4$.
\end{proof}

So, by Properties (1),(2), and (3) for any $\lambda \in
(D^2)^\perp\backslash (D^3)^\perp$ we get the flag
\begin{equation}
\label{flag}
\ldots\subseteq\mJ_{(i)}(\lambda)\subseteq\ldots\subseteq\mJ_{(1)}(\lambda)\subset
\mJ(\lambda)\subset\mJ^{(1)}(\lambda)\subseteq\ldots\subseteq
\mJ^{(i)}(\lambda)\subseteq\ldots
\end{equation}
in $T_\lambda (D^2)^\perp$.
 %define the additional
%%additional monotonically nonincreasing
%sequence of distributions $\{{\mathcal J}_{(i)}\}_{i\geq 0}$ on
%$(D^2)^\perp\backslash (D^3)^\perp$ as  follows:
%
%\begin{equation}
%\label{contr} {\mathcal J}_{(i)}(\lambda)= \{v\in T_\lambda
%\bigl((D^2)^\perp\bigr): \sigma (v, w)=0\,\, \forall w\in {\mathcal
%J}^{(i)}(\lambda)\},
%\end{equation}
%%In other words, the subspace ${\mathcal J}_{(i)}(\lambda)$ is the
%%skew-symmetric complement of ${\mathcal J}^{(i)}(\lambda)$ (w.r.t.
%%the form $\sigma$).
The dynamics of these flags
%$(D^2)^\perp$)
along any abnormal extremal
%of $D$
%characteristic curve (of  $W_D$)
%$(D^2)^\perp$)
defines certain
%unparameterized
curve of flags of isotropic and coisotropic subspaces in a linear
symplectic space.

More precisely, let $\gamma$ be a segment of abnormal extremal of
$D$ and $O_\gamma$ be a neighborhood of $\gamma$ in $(D^2)^\perp$
such that the factor
%\begin{equation}
%\label{Ndef}
$N=O_\gamma /(\text {\emph{the characteristic one-foliation}})$
%\end{equation}
 is a
well defined smooth manifold. The quotient manifold $N$ is a
symplectic manifold endowed with the symplectic structure
$\bar\sigma$ induced by $\sigma |_{(D^2)^\perp}$. Let $\phi
:O_\gamma\to N$ be the canonical projection on the factor.
%Its dimension is equal to $2(n-2)$.
%The quotient manifold $N$ is
%a symplectic manifold endowed with a symplectic structure
%$\bar\sigma$ induced by $\sigma |_{(D^2)^\perp}$.
%Let $\phi \colon O_\gamma\to N$ be the canonical projection on the
%factor.
%It is easy to check that $\phi_*\bigl({\mathcal
%J}(\lambda)\bigr)$ is a Lagrangian subspace of the
%symplectic space $T_{\gamma}N$ for all $\lambda\in\gamma$.
 %Let $L(T_\gamma
%N)$ be the Lagrangian Grassmannian of the symplectic space
%$T_\gamma N$, i.e., $L(T_\gamma N)=\{\Lambda\subset
%T_\gamma N:\Lambda^\angle=\Lambda\}$, where
%$\Lambda^\angle$ is the skew-symmetric complement of the
%subspace $\Lambda$, $\Lambda^\angle=\{v\in T_\gamma
%N:\bar\sigma(v,\Lambda)=0\}$.
%\begin{defin}
%\label{jacdef}
For each $\ge0$ we can define the following curves of subspaces in $T_\gamma N$:
%the mapping
%% \emph{ Jacobi curve of $\gamma$} is the mapping
%$J_\gamma\colon\gamma\mapsto G_{n-2}(T_\gamma N)$ by
%%be the mapping such that
\begin{equation}
\label{jacurve}
\lambda\mapsto
%\widetilde J
%%_\gamma
%%^{(i)}
%(\lambda)
%%\stackrel{def}
%{=} \phi_*\bigl({\mathcal J}
%%^{(i)}
%(\lambda) \bigr), \quad
%\widetilde J
%_\gamma
%^{(i)}(\lambda)
%\stackrel{def}
%{=}
\phi_*\bigl({\mathcal J}^{(i)}(\lambda) \bigr), \quad
\lambda\mapsto
%\widetilde J_{
%\gamma
%(i)}(\lambda)
%\stackrel{def}
%{=}
\phi_*\bigl({\mathcal J}_{(i)}(\lambda) \bigr),\quad
 \text {for all }
 %i\in \mathbb N,
\lambda\in\gamma.
\end{equation}
%For shortness we also set $\widetilde J(\lambda)=\widetilde
%J^0(\lambda)$ ($=\widetilde J_0(\lambda)$). First note that
These curves describe the dynamics of the corresponding subspaces of the flag
\eqref{flag} w.r.t. the characteristic $1$-foliation along the abnormal extremal
$\gamma$.

Note that there exists a straight line, which is common to all subspaces
appearing in \eqref{jacurve} for any $\lambda\in\gamma$. So, it is more
convenient to get rid of it by a factorization. Indeed, let $e$ be the Euler
field on $T^*M$, i.e., the infinitesimal generator of homotheties on the fibers
of $T^*M$. Since a transformation of $T^*M$, which is a homothety on each fiber
with the same homothety coefficient, sends abnormal extremals to abnormal
extremals, we see that the vector $\bar e= \phi_*e(\lambda)$ is the same for any
$\lambda\in \gamma$ and lies in any subspace appearing in \eqref{jacurve}. Let
\begin{equation}
\label{jacurve1} J^{(i)}(\lambda)=\phi_*\bigl({\mathcal
J}^{(i)}(\lambda)\bigr)/\{\mathbb R \bar e\},\quad
J_{(i)}(\lambda)=\phi_*\bigl({\mathcal
J}_{(i)}(\lambda)\bigr)/\{\mathbb R \bar e\}.
\end{equation}
For simplicity we also set $J(\lambda)=J^{(0)}(\lambda)$ ($=
J_{(0)}(\lambda)$). It is clear that all subspaces appearing in
\eqref{jacurve1} belong to the space
\begin{equation}
\label{spaceW}
 W=\{v\in T_\gamma N: \bar\sigma(v,\bar
e)=0\}/\{\mathbb R \bar e\}.
\end{equation}
The space $W$ is endowed with the natural symplectic structure
induced by $\bar\sigma$, which for simplicity will be denoted also
by $\bar\sigma$. Also $\dim W=2(n-3)$. Rewriting Properties (1)-(4)
of Proposition \ref{proper} and relation (\ref{contr}) in terms of
subspaces $J^{(i)}(\lambda)$ and $J_{(i)}(\lambda)$, we get

\begin{prop}
\label{proper1} ~

\begin{enumerate}
\item The subspace $J(\lambda)$ is  Lagrangian
subspace of $W$ and  the subspace $J_{(i)}(\lambda)$ is the
skew-symmetric complement of $J^{(i)}(\lambda)$ in $W$ for any
$\lambda\in\gamma$;

\item $J^{(i-1)}(\lambda)\subseteq  J^{(i)}(\lambda)$,
 $J_{(i)}(\lambda)\subseteq J_{(i-1)}(\lambda)$;

\item ${\rm dim}\,J^{(1)}(\lambda)-{\rm dim}\, J(\lambda) =1$,  ${\rm dim}\,
J(\lambda)-{\rm dim}\,  J_{(1)}(\lambda) =1$;

\item ${\rm dim}\, J^{(i)}(\lambda)-{\rm dim}\,  J^{(i-1)}(\lambda)\leq 1$, ${\rm dim}\,
J_{(i-1)}(\lambda)-{\rm dim}\, J_{(i)}(\lambda)\leq 1$ for $i\ge2$;

%\item ${\rm dim}\, {\mathcal J}^{(i)}(\lambda)\leq 2n-4$.
%%\item$\{\mathcal C\}, e\}
\end{enumerate}
\end{prop}

%, where ${\mathcal J}(\lambda)$ is as above. Actually, the
%symplectic form $\sigma$ of $T^*M$ induces naturally the symplectic
%form $\bar \sigma$ on $T_\gamma N$ and
Note also that  by Property (1) of Proposition \ref{proper} for all
$i\in \mathbb N$ the subspaces $J^{(i)}(\lambda)$ are coisotropic
and the subspaces $J_{(i)}(\lambda)$ are isotropic in $T_\gamma N$.
%So, segment $\gamma$ of abnormal extremal $\gamma$ one can construct
The curve
\begin{equation}
\label{flag1} \lambda\mapsto \left\{ \ldots\subseteq
J_{(i)}(\lambda)\subseteq\ldots\subseteq J_{(1)}(\lambda)\subset
J(\lambda)\subset J^{(1)}(\lambda)\subseteq\ldots\subseteq
J^{(i)}(\lambda)\subseteq\ldots \right\}, \quad \lambda \in \gamma,
\end{equation}
of flags of isotropic and coisotropic subspaces in a linear
symplectic space $W$ will be called the \emph{curve of flags
associated with the segment $\gamma$ of an abnormal extremal}.

Clearly, any symplectic invariant of such curve automatically produces an
invariant of the distribution $D$ itself. Moreover, it turns out that under
certain generic assumptions one can construct the canonical frames of the
distribution $D$ from the study of differential geometry of such curves.

\begin{remark}
\label{reclagr} {\rm As a matter of fact the whole curve of flags
\eqref{flag1} can be recovered from the curve $\lambda\mapsto
J(\lambda)$ of Lagrangian subspaces of $W$. It is clear by item 1 of
Proposition \ref{proper1} that it is sufficient to show how to
recover the subspaces $J^{(i)}(\lambda)$. This can be done by making
an appropriate differentiation (in a similar manner as all subspaces
${\mathcal J}^{(i)}(\lambda)$ and ${\mathcal J}_{(i)}(\lambda)$ are
obtained from ${\mathcal J}(\lambda)$). More precisely, let
$\Gamma(J)$ be the set of all smooth mappings $\ell$ from $\gamma$
to the ambient symplectic space $W$ (see \eqref{spaceW}) such that
$\ell(\lambda)\in J(\lambda)$ for all $\lambda\in\gamma$. In other
words, $\Gamma(J)$ is the space of all smooth sections of the vector
bundle over $\gamma$ having the subspace $J(\lambda)$ as the fiber
over a point $\lambda\in\gamma$. If $\varphi:\gamma\mapsto\mathbb R$
is a parameterization of $\gamma$, $\varphi(\lambda)=0$ and
$\psi=\varphi^{-1}$, then
\begin{equation}
\label{Jispan} J^{(i)}(\lambda)={\rm span}\Bigl\{
\frac{d^j}{dt^j}\ell\bigl(\psi(t)\bigr)|_{t=0}: \ell\in \Gamma(J),\,
% \text
% {is a smooth section of}\mathfrak G,
0\leq j\leq i\Bigr\}.
\end{equation}
The curve $\lambda\mapsto J(\lambda)$ is called \emph {Jacobi curve
of the abnormal extremal $\gamma$}. The reason to call this curve
Jacobi curve is that it can be considered as the generalization of
spaces of "Jacobi fields" along Riemannian geodesics: in terms of
this curve one can describe some optimality properties (so-called
rigidity) of the corresponding abnormal extremal (\cite{agrachev} or
\cite{zel})}.$\Box$
\end{remark}

\begin{remark}
\label{otherrep} \rm{Actually, one can describe the subspaces
$J_{(i)}(\lambda)$, where $i\geq 1$, without using the symplectic
structure on $W$. For this, by analogy with above, let
$\Gamma(J_{(i)})$, $i\geq 0$,  be the set of all smooth mappings
$\ell$ from $\gamma$ to the ambient symplectic space $W$  such that
$\ell(\lambda)\in J_{(i)}(\lambda)$ for all $\lambda\in\gamma$. Let
$\varphi:\gamma\mapsto\mathbb R$ be a parameterization of $\gamma$,
$\varphi(\lambda)=0$ and $\psi=\varphi^{-1}$. Then it is easy to
show that for any $i\geq 1$
\begin{equation}
\label{contr1} J_{(i)}(\lambda)=
%{\mathcal D}^{(i-1)}\Lambda (\tau)+
\left\{v\in  J_{(i-1)}(\lambda):
\begin{array}{l}
\exists\,\ell\in\Gamma(J_{(i-1)})\,\,
%\text {a vector field}\,\,
%{\mathcal V}\subset {\mathcal J}_{(i-1)}\,\,
\text {with}\,\, \ell(\lambda)=v\\
%{\mathcal V}(\lambda)=v
\text {such that}\,\,
%\bigr[{\mathcal C}, {\mathcal
%V}\bigl](\lambda)\in {\mathcal J}_{(i-1)}(\lambda),
\frac{d}{dt}\ell\bigl(\psi(t)\bigr)|_{t=0}\in
J_{(i-1)}(\lambda)
\end{array}\right\}
\end{equation}
The last formula allows to construct $J_{(i)}$ recursively, starting
from $J_{(0)}=J$. Finally it is easy to show  that identity
\eqref{contr1} remains true if one replaces the quantor $\exists$ by
$\forall$. }$\Box$
\end{remark}
\section {The class of a rank 2 distribution}
\setcounter{equation}{0}

 Let us describe precisely the generic assumptions on a germ of rank $2$ distribution
 necessary for constructing the canonical frames for them.

First for any point $q\in M$ denote by $(D^l)^\perp(q)
=(D^l)^\perp\cap T_q^*M$ the fiber of $(D^l)^\perp$. Let us define
the following
%two
integer-valued function on
$(D^2)^\perp\backslash(D^3)^\perp$
%and on $M$ respectively
:
\begin{equation*}
\begin{split}
%&~&\label{K}
~& \nu(\lambda)=\min\{i\in {\mathbb N}: {\mathcal
J}^{(i+1)}(\lambda)={\mathcal J}^{(i)}(\lambda)\}, \\
%and
%\end{equation}
%\begin{defin}
%&~&\label{class }
~&m(q)=\max\{\nu(\lambda):\lambda\in
(D^2)^\perp(q)\backslash(D^3)^\perp(q)\}.
\end{split}
\end{equation*}
From  Properties (3), (4), (5) of Proposition \ref{proper}, and the
fact that $\dim \mJ = n-1$ it follows that $1\leq \nu(\lambda)\leq
n-3$. Furthermore, from \eqref{foli25} and the definition
\eqref{Ji1} of $\mJ_{(i)}$ it follows that the set $\{\lambda\in
(D^2)^\perp(q): \nu(\lambda)=m(q)\}$ is nonempty open in the Zariski
topology of the fiber $(D^2)^\perp(q)$. Besides, it is easy to show
that the integer-valued functions $\nu(\cdot)$ and $m(\cdot)$ are
lower semicontinuous. Hence they are locally constant on the open
and dense subset of $(D^2)^\perp\backslash (D^3)^\perp$ and $M$
correspondingly
 and attain their maximum values on the open sets there.
% Obviously, the set of points, where the
%distribution $D$ has maximal class $n-3$, is open.
Moreover,

\begin{prop}
\label{nonhol}
 Germs of $(2,n)$-distributions of the maximal class $n-3$
 are generic.
\end{prop}
This Proposition was proved in \cite[Proposition 3.4]{zelvar}. Let
us outline the proof. The distribution has maximal class at a point
$q_0$ if and only if its jet of sufficiently high order belongs to
the Zariski open set of the jet space of this order. Therefore in
order to prove Proposition \ref{nonhol} it is sufficient to give
just one example of a germ of $(2,n)$-distributions of the maximal
class $n-3$. As such example one can take the germ at $0$ of the
distribution, associated with underdetermined ODE \eqref{under},
where $F=\frac{1}{2}p_{n-3}^2$.

In the present paper we treat the germs of $(2,n)$ distributions of the maximal
class $n-3$. In the cases $n=5$ and $n=6$ any rank 2 distribution has maximal
class if and only if it has maximal possible small growth vector, namely,
$(2,3,5)$ in the case $n=5$ and $(2,3,5,6)$ in the case $n=6$. It can be obtained
by direct computations (see Propositions 3.5 and 3.6 of \cite{zelvar}
respectively). Starting with $n=7$, distributions with different small growth
vectors may have the maximal class.

Regarding completely nonholonomic rank 2 distributions of
non-maximal class, it is easy to describe all distributions of
minimal class $1$. From \cite[Remark~3.4]{zelvar} it follows that a
rank~2 distribution $D$ has the smallest possible class $1$ at a
point $q$ if and only if $\dim\, D^3(q)=4$. Moreover, it is easy to
see that in this case $D$ is either the Goursat distribution or, by
the factorization of the ambient manifold by the characteristics of
$D^2$ (or series of such factorizations), one can get a new
distribution $\widetilde D$, satisfying $\dim \widetilde D^3=5$ (or,
equivalently, having the class greater than $1$). See the last
section of the paper for more details.

On the other hand, \textbf{\emph{we have no examples of completely nonholonomic
rank $2$ distributions of constant class $2\leq m\leq n-4$ and our conjecture is
that there are no such distributions}}. We succeeded to prove this conjecture for
$m=2,3,4$ by direct computation, which means, in particular, that any
counter-example to our conjecture, if it exists, should live on at least a
$9$-dimensional manifold.

\begin{remark}
\label{projCr} {\rm By above, for $(2,n)$-distributions of maximal class with
$n>4$ it is necessary that $\dim\, D^3=5$. Hence, each such distribution
satisfies relation (\ref{projC}).} $\Box$
\end{remark}

From now on $D$ is a $(2,n)$-distribution of maximal class $m=n-3$. Let us study
curves of flags associated with its abnormal extremals in more detail. Let
\begin{equation}
\label{Jconstw2} {\mathcal R}_D=\{\lambda\in (D^2)^\perp \backslash
(D^3)^\perp: \nu(\lambda)=n-3 \}, \quad  {\mathcal R}_D(q)={\mathcal
R}_D\cap T_q^*M.
\end{equation}
As was already mentioned the set ${\mathcal R}_D$ is open dense
subset of $(D^2)^\perp \backslash (D^3)^\perp$ and the set
${\mathcal R}_D(q)$ is a nonempty open set in Zariski topology on
the linear space $(D^2)^\perp(q)$.
The following Proposition follows easily from Proposition \ref{proper}, formula \eqref{Jispan}, and
Remark \ref{otherrep}
%(see
%again~\cite[Proposition~3.4]{zelvar}).

\begin{prop}
\label{projcurve}
%Suppose that a point $\bar\lambda\in
%(D^2)^\perp\backslash (D^2)^\perp$ satisfies $\nu(\bar\lambda)=n-3$.
Let $\gamma$ be a segment of abnormal extremal such that
$\gamma\subset {\mathcal R}_D$.
%passing through a
%point $\bar\lambda$ with $\nu(\bar\lambda)=n-3$.
Then
%the germ  of
the associated curve of flags is the curve of complete flags in $W$,
i.e., it has the form
\begin{equation}
\label{flag2}
\begin{split}
~&\lambda\mapsto \{ 0= J_{(n-3)}(\lambda)\subset
J_{(n-4)}(\lambda)\subset\ldots\subset J_{(1)}(\lambda)\subset
J(\lambda)\subset\\ ~& J^{(1)}(\lambda)\subset\ldots\subset
J^{(n-4)}(\lambda)\subset J^{(n-3)}(\lambda)=W\}, \quad \lambda \in
\gamma,
\end{split}
\end{equation}
where $\dim J_{(i)}(\lambda)=n-3-i$ and $\dim
J^{(i)}(\lambda)=n-3+i$.
%(for $\lambda$ close to $\bar\lambda)$.

Moreover, if $\lambda\mapsto\ell(\lambda)$ is a smooth curve of
vectors such that $J_{(n-4)}(\lambda)=\mathbb R\ell(\lambda)$ ,
$\varphi\colon\gamma\mapsto\mathbb R$ is a parameterization of
$\gamma$,
%$\varphi(\bar\lambda)=0$
and $\psi=\varphi^{-1}$, then for $i=0,\ldots n-3$ and any $t\in
\varphi(\gamma)$
%sufficiently close to $0$
\begin{equation}
\label{diffone}
\begin{split}
~& J_{(i)}\bigl(\psi(t)\bigr)={\rm
span}\bigl\{\frac{d^j}{dt^j}\ell\bigl(\psi(t)\bigr):0\leq
j\leq n-4-i\bigr\}\\
~& J^{(i)}\bigl(\psi(t)\bigr)={\rm
span}\bigl\{\frac{d^j}{dt^j}\ell\bigl(\psi(t)\bigr):0\leq j\leq
n-4+i\bigr\}.
\end{split}
\end{equation}
\end{prop}

In other words,
%the previous proposition shows that
the curve of flags,
associated with the abnormal extremal $\gamma\subset \mathcal R_D$
can be recovered by differentiation from the curve of their
one-dimensional subspaces $\lambda\mapsto J_{(n-4)}(\lambda)$, i.e.,
from the curve in the projective space $\mathbb P W$ of the
$2m$-dimensional symplectic space $W$  ($m=n-3$). Moreover, the
curve $\lambda\mapsto J_{(n-4)}(\lambda)$ is not arbitrary curve in
$\mathbb P W$ but a curve, which can be completed by the appropriate
number of differentiation to the curve of Lagrangian subspaces of
$\gamma$.

The differential geometry of curves in projective spaces is the
classical subject, essentially completed already in 1905 by
E.J.~Wilczynski (\cite{wil}). In particular, it is well known that
these curves  are endowed with the canonical projective structure,
i.e., we have a distinguished set of parameterizations (called
projective) such that the transition function from one such
parameterization to another is a M\"{o}bius transformation. Let us
demonstrate how to construct it for the curve $\lambda\mapsto
J_{(n-4)}(\lambda)$, $\lambda\in\gamma$.

As before, let $\Gamma(J_{(n-4)}(\lambda))$ be the space of all
smooth mappings $\ell$ from $\gamma$ to the ambient symplectic space
$W$ such that $\ell(\lambda)\in J_{(n-4)}(\lambda)$ for all
$\lambda\in\gamma$. The elements of $\Gamma(J_{(n-4)}(\lambda))$
will be called \emph{sections of the curve $\lambda\mapsto
J_{(n-4)}(\lambda)$}. Take some parameterization
$\varphi\colon\gamma\mapsto {\mathbb R}$ of $\gamma$ and let
$\psi=\varphi^{-1}$. By Proposition \ref{projcurve} for any section
$\ell$ one has relation
\begin{equation}
\label{spanall} {\rm span}\bigl\{\frac{d^j}{dt^j}\ell\bigl(\psi(t)\bigr)\mid
0\leq j\leq 2m-1\bigr\}=W.
\end{equation} It is well known that there exists the unique, up to the
multiplication on a nonzero constant, section
%$t\mapsto
$E_\varphi$ such that
%$(t)$,
%corresponding to the parameter $t$
%rization
%$E$if $\tau$ is another parameter,
%$t\mapsto E(t)$
\begin{equation}
\label{last0}
\frac{d^{2m}}{dt^{2m}}E_\varphi\bigl(\psi(t)\bigr)=\sum_{i=0}^{2m-2}B_i^\varphi(t)\frac{d^{i}}{dt^{i}}
E_\varphi\bigl(\psi(t)\bigr),
\end{equation}
i.e. the coefficient of the term
$\frac{d^{2m-1}}{dt^{2m-1}}E_\varphi\bigl(\psi(t)\bigr)$ in the linear
decomposition of $\frac{d^{2m}}{dt^{2m}}E_\varphi\bigl(\psi(t)\bigr)$ w.r.t. the
basis $\bigl\{\frac{d^i}{dt^i}\ell\bigl(\psi(t)\bigr):0\leq i\leq 2m-1\bigr\} $
vanishes.

Further, let $\varphi_1$ be another parameter and
$\upsilon=\varphi\circ\varphi_1^{-1}$. Then it is not hard to show
that the coefficients and $B_{2m-2}^\varphi$ and
$B_{2m-2}^{\varphi_1}$ in the decomposition (\ref{last0}),
corresponding to parameterizations $\varphi$ and $\varphi_1$, are
related as follows:
%the canonical representative of
%$\zeta$ w.r.t. the parameter $\tau$, then
\begin{equation}
 \label{rhorep}
 \widetilde
 B_{2m-2}^{\varphi_1}(\tau)=\upsilon'(\tau)^2B_{2m-2}^{\varphi}(\upsilon(\tau))-\frac{m(4m^2-1)}{3}
\mathbb{S}(\upsilon)(\tau),
\end{equation}
 where $\mathbb{S}(\upsilon)$ is a
Schwarzian derivative of $\upsilon$,
%\begin{equation}
%\label{sch}
$\mathbb S(\upsilon)=
%\frac{1}{2}\frac{\varphi^{(3)}}{\varphi'}-
%\frac{3}{4}\Bigl (\frac{\varphi''}{\varphi'}\Bigr)^2
\frac {d}{dt}\Bigl(\frac {\upsilon''}{2\,\upsilon'}\Bigr)
-\Bigl(\frac{\upsilon''}{2\,\upsilon'}\Bigr)^2$.

From the last formula and the fact that ${\mathbb S}\upsilon\equiv
0$ if and only if the function $\upsilon$ is M\"{o}bius it follows
that \emph{the set of all parameterizations $\varphi$ of $\gamma$
such that
\begin{equation}
\label{A2m-2} B_{2m-2}^{\varphi}\equiv 0
\end{equation}
defines the canonical projective structure on $\gamma$}. Such
parameterizations are called the \emph {projective parameterizations
of the abnormal extremal $\gamma$}.

\begin{remark}
\label{crossrat} {\rm Another description of the canonical
projective structure on an abnormal extremal $\gamma$ can be
obtained by working with the Jacobi curve $\lambda\mapsto
J(\lambda)$, using the notion of the cross-ratio of four point in
Lagrangian Grassmannian (\cite{jac1},\cite{zelvar}). This approach
allows to construct the canonical projective structures in much more
general situations.}
\end{remark}

Further, since in our case the $m$-dimensional subspaces
$J\bigl(\psi(t)\bigr)$ are Lagrangian, it is easy to show that the
condition \eqref{last0} for the section $E_\varphi(t)$ is equivalent
to the following one
%the sections $E_\varphi$ (defined up to a m
%have the following representative $t\mapsto \varepsilon(t)$ is
%canonical if and only if
\begin{equation}
\label{symp} \bar\sigma
\left(\frac{d^{m}}{dt^{m}}E_\varphi\bigl(\psi(t)\bigr),\frac{d^{m-1}}{dt^{m-1}}
E_\varphi\bigl(\psi(t)\bigr)\right)\equiv {C},\,\quad C\in\mathbb
R\backslash\{0\}.
\end{equation}
Therefore in our case we can ``kill'' the freedom of the multiplication on a
nonzero constant in the definition of section $E_\varphi$ by setting

\begin{equation}
\label{norme1}
\left|\bar\sigma\left(\frac{d^{m}}{dt^{m}}E_\varphi\bigl(\psi(t)\bigr),\frac{d^{m-1}}{dt^{m-1}}
E_\varphi\bigl(\psi(t)\bigr)\right)\right|\equiv 1.
\end{equation}
(see also \cite{praga}). There are exactly two  sections of the
curve $\lambda\mapsto J_{(n-4)}(\lambda)$ with the parameterization
$\varphi$, satisfying (\ref{norme1}), and they are obtained one from
another by the multiplication on $-1$. These sections are called
\emph{the canonical sections of the curve $\lambda_\mapsto
J_{(n-4)}(\lambda)$ w.r.t. the  parameterization $\varphi$} and both
of them  will be denoted in the sequel by $E_\varphi$. From these
sections one can obtain the moving frame
$\bigl\{\frac{d^i}{dt^i}E_\varphi\bigl(\psi(t)\bigr):0\leq i\leq
2m-1\bigr\}$ on $W$, which is defined again up to a multiplication
by $-1$ and it will be called \emph{the canonical moving frame of
the curve $\lambda\mapsto J_{(n-4)}(\lambda)$ w.r.t. the
parameterization~$\varphi$}.

Finally, it can be shown easily that in the case when $m$-dimensional subspaces
$J\bigl(\psi(t)\bigr)$ are Lagrangian, from \eqref{A2m-2} it follows that
\begin{equation}
\label{A2m-3} B_{2m-3}^{\varphi}\equiv 0
\end{equation}

\section{The canonical frame}
\setcounter{equation}{0}

Now we are ready to describe the manifold, on which the
canonical frame for $(2,n)$-distribution of maximal class,
$n>5$, can be constructed. Given $\lambda\in {\mathcal
R}_D$ denote by ${\mathfrak P}_\lambda$ the set of all
projective parameterizations $\varphi:\gamma\mapsto\mathbb
R$ on the characteristic curve $\gamma$ , passing through
$\lambda$, such that $\varphi(\lambda)=0$. Denote
%\begin{equation*}
$$\Sigma_D=\{(\lambda, \varphi):\lambda\in {\mathcal R}_D,
\varphi\in {\mathfrak P}_\lambda\}.$$ Actually, $\Sigma_D$
is a principal bundle over ${\mathcal R}_D$ with the
structural group of all M\"{o}bius transformations,
preserving $0$ and $\dim \, \Sigma_D=2n-1$.
%Actually,
%${\mathcal H}$ is isomorphic to the subgroup of lower
%triangular matrices of $SL(2,\mathbb R)$. Our main result
%is the following

\begin{theor}
\label{frtheor}
%\label{main}
For any $(2,n)$-distribution, $n>5$, of maximal class
there exist two canonical frames on the corresponding
$(2n-1)$-dimensional manifold $\Sigma_D$, obtained one from
another by a reflection. The group of symmetries of such
distributions is at most $(2n-1)$-dimensional.
\end{theor}
%, where the action of $\mathfrak{gl}(2,\mathbb
%R)$ on a subspace of ${\mathfrak n}_{2n-5}$, which is
%complement to its centre, is isomorphic to the standard
%irreducible representation of $\mathfrak{gl}(2,\mathbb R)$
%on ${\rm Symm}^{n-3}$.
%\end{theor}
%\newpage

{\bf Proof.} Define the following two
fiber-preserving flows on $\Sigma_D$:
\begin{equation}
\label{flows} F_{1,s}(\lambda,\vf)=(\lambda,
e^{2s}\varphi),\quad
F_{2,s}(\lambda,\varphi)=\left(\lambda,
\cfrac{\varphi}{-s\varphi+1}\right),\quad
\lambda\in{\mathcal R}_D, \varphi\in {\mathfrak P}_\lambda.
\end{equation}

 Further, let $\delta_s$ be the flow of homotheties on
 the fibers of $T^*M$:
\begin{equation}
\label{homoth} \delta_s(p,q)=(e^sp,q),\quad q\in M,\,\,
p\in T_q^*M
\end{equation}
(actually the Euler field $e$ generates this flow). The
following flow
\begin{equation}
\label{flow0}
F_{0,s}(\lambda,\varphi)=\bigl(\delta_{2s}(\lambda),
\varphi\circ\delta_{2s}^{-1}\bigr)
\end{equation}
is well-defined on $\Sigma_D$ (here we use that $\delta_s$
preserves the characteristic $1$-foliation).
%\end{equation}
For any $0\leq i\leq 2$ let $g_i$ be the vector field on
$\Sigma_D$, generating the flow $F_{i,s}$. Note that $g_1$ and $g_2$ are just
fundamental vector fields on $\Sigma_D$ defined by the structure of the
principle fiber bundle on $\Sigma_D$.

%Let $g_0$ be the vector
%field on $\Sigma$ such that
%$g_0(u)=\frac{\partial}{\partial s}F_{0,s} (u)|_{s=1}$ for
%any $u\in \Sigma$.
Besides, the characteristic $1$-foliation on $(D^2)^\perp$ can be
lifted to the \emph{parameterized} $1$-foliation on $\Sigma_D$,
which gives one more canonical vector field on $\Sigma_D$. Indeed,
let $u=(\lambda,\varphi)\in \Sigma_D$ and $\gamma$ be the
characteristic curve, passing through $\lambda$ (so, $\vf$ maps
$\gamma$ to $\mathbb R$). Then the mapping
$$\Upsilon_u(t)=\bigl(\vf^{-1}(t), \vf(\cdot)-t\bigr)$$ defines the
parameterized curve on $\Sigma_D$, the lift of $\gamma$ to
$\Sigma_D$, and $\Upsilon_u(0)=u$. The additional canonical vector
field $h$ on $\Sigma_D$ is defined by
\begin{equation}
\label{flowh} h(u)=\frac{d}{dt}\Upsilon_u(t)|_{t=0}.
\end{equation}
It can be shown easily that
\begin{equation}
\label{gl2}
[g_1,g_2]=2g_2,\,\, [g_1, h]=-2h,\,\, [g_2,h]=g_1,\,\,
[g_0, h]=0,\,\,[g_0,g_i]=0.
\end{equation}
%and $g_0$ commutes with all constructed fields.
Therefore
the linear span (over ${\mathbb R}$) of the vector fields
$g_0$, $g_1$, $g_2$, and $h$ is endowed with a structure of the
Lie algebra isomorphic to
$\mathfrak{gl}(2,\mathbb{R})$.

Now we will construct one more canonical, up to
the sign, vector field on $\Sigma_D$. %For this let
%\begin{equation}
%\label{contr} {\mathcal J}_{(i)}(\lambda)= \{v\in
%T_\lambda \bigl((D^2)^\perp\bigr): \sigma (v, w)=0\,\,
%\forall w\in {\mathcal J}^{(i)}\}.
%\end{equation}
%First from (\ref{simpJ}) it follows that ${\mathcal
%J}_{(0)}={\mathcal J}$.
First from (\ref{J1}) it is easy to get
\begin{equation}
\label{J_1} \mJ_{(1)}(\lambda)=
T_\lambda\Bigl((D^2)^\perp\bigl(\pi(\lambda)\bigr)\Bigr)\oplus
{\mathcal C}
\end{equation}
Here $T_\lambda\Bigl((D^2)^\perp\bigl(\pi(\lambda)\Bigr)$ is the
tangent space to the fiber $(D^2)^\perp\bigl(\pi(\lambda)\bigr)$ at
the point $\lambda$ and it is actually equal to $\{v\in T_\lambda
(D^2)^\perp,\pi_*v=0\}$, the vertical subspace of $T_\lambda
(D^2)^\perp$.
%Besides, directly from (\ref{contr}) and definitions of the characteristic foliation ${\mathcal C}$ and the Euler field $e$ it follows that
%\begin{equation}
%\label{charein}
%\forall i\in {\mathbb N},\quad {\rm span} \{C,e\}\subset J_{(i)}.
%\end{equation}
Let
%$\sigma|_{\mathcal J}=0$.
\begin{equation}
\label{vertcomp} V_i(\lambda)={\mathcal
J}_{(i)}(\lambda)\cap
T_\lambda\Bigl((D^2)^\perp\bigl(\pi(\lambda)\bigr)\Bigr).
%\end{equation}
\end{equation}
Since
%${\mathcal
%J}^{(i)}\subseteq\mJ^{(i+1)}$, then
${\mathcal J}_{(i+1)}\subseteq\mJ_{(i)}$,
%. Therefore by
identity (\ref{J_1}) yields
%and (\ref{charein})
\begin{equation}
\label{split}  {\mathcal J}_{(i)}=V_i\oplus {\mathcal C}
\quad\forall i\geq 1.
\end{equation}
Note also that %if by analogy with above, $\Gamma(V_{i})$ is the set of
%all smooth mappings $\ell$ from $\gamma$ to the ambient symplectic
%space $W$ such that $\ell(\lambda)\in V_{i}(\lambda)$ for all
%$\lambda\in\gamma$, then f
from Remark \ref{otherrep} and formula \eqref{split} it follows that

\begin{equation}
\label{contr2} V_i(\lambda)=
%{\mathcal D}^{(i-1)}\Lambda (\tau)+
\left\{v\in V_{i-1}(\lambda):
\begin{array}{l}\exists\,\, \text {a vector field}\,\,
{\mathcal V}\in V_{i-1} \,\, \text {with}\\
{\mathcal V}(\lambda)=v \,\,\text {such that }\,\, \bigr[{\mathcal
C}, {\mathcal V}\bigl](\lambda)\in {\mathcal J}_{(i-1)}(\lambda)
\end{array}\right\}.
\end{equation}
Furthermore, identity \eqref{contr2} remains true if one
replaces the quantor $\exists$ by $\forall$. In particular, from the
last identity and Proposition \ref{projcurve} it follows that for
any $\lambda\in{\mathcal R}_D$ the distributions ${\mathcal
J}_{(i)}$ and $V_i$ satisfy
\begin{equation}
\label{hV} {\mathcal J}_{(i)}(\lambda)=[{\mathcal
C},\mJ_{(i+1)}],\quad V_{i}(\lambda)=[{\mathcal
C},V_{i+1}], \quad 0\leq i< n-4
%{\mathcal D}^{(i-1)}\Lambda (\tau)+
%\left\{v\in {\mathcal J}_{(i-1)}(\lambda):
%\begin{array}{l}\exists\,\, \text {a vector field}\,\, {\mathcal V}\subset {\mathcal
%J}_{(i-1)}\,\, \text {such that}\\ {\mathcal V}(\lambda)=v
%\,\,{\rm and }\,\, \bigr[{\mathcal C}, {\mathcal
%V}\bigl](\lambda)\in {\mathcal J}_{(i-1)}(\lambda)
%\end{array}\right\}
\end{equation}
Take vector fields $E\in V_{n-4}$ and $H\in {\mathcal C}$ without
stationary points and suppose that $E$ is not collinear to the Euler
field $e$. Then by our construction of subspaces $J_{(i)}$ and
$J^{(i)}$, Proposition \ref{projcurve} and relation \eqref{split} it
follows that on $\mathcal R_D$
\begin{equation}
\label{basVJ}
\begin{aligned}
~&V_{n-4}=%{\rm span}\,
\langle
%\{
e, E%\}
\rangle \, \quad \mJ_{(i)}=
%{\rm span}\left\{
\langle H, e, E, \{({\rm ad} H)^j E\}_{j=1}^{n-4-i} \rangle
%\right\},
\quad
\,\, 0\leq i\leq n-4,\\
~&\mJ^{(i)}=
%{\rm span}\left\{
\langle H, e, E, \{({\rm ad} H)^j E\}_{j=1}^{n-4+i} \rangle
% \}
\,\, 0\leq i\leq n-3.
\end{aligned}
\end{equation}

%In particular, if $\lambda\in\mathcal R_D(\lambda)$ , then from the
%last identity, Proposition \ref{projcurve} and relation
%\eqref{split} we will get immediately that $\dim V_{n-4}=2$.
%Besides, by our constructions the Euler field $e$ is tangent to rank
%2 distribution $V_n-4$.
Now let $\gamma$ be the abnormal extremal, passing through
$\lambda\in\mathcal R_D(\lambda)$. As before, let also $\phi$ be the
canonical projection from a sufficiently small neighborhood
$O_\gamma$ of $\gamma$ to the factor $O_\gamma /(\text {\emph{the
characteristic one-foliation}})$, $\bar e=\phi_* e(\lambda)$, and
$E_\varphi(\lambda)$ be one of the two canonical sections of the
curve $\lambda_\mapsto J_{(n-4)}(\lambda)$ w.r.t. the
parameterization $\varphi$.  Then, obviously, there exists a unique
affine line ${\rm Aff}_\varphi(\lambda)$ in the plane
$V_{n-4}(\lambda)$ such that
\begin{equation}
\label{aff} \phi_*\bigl({\rm Aff}_\varphi(\lambda)\bigr)/\{\mathbb R
 e\}=E_\varphi(\lambda).
 \end{equation}
Clearly, the affine line ${\rm Aff}_\varphi(\lambda)$ is parallel to
the vector $e(\lambda)$, but does not pass through the origin of the
linear space $V_{n-4}(\lambda)$.

Further, denote by $\Pi:\Sigma_D\mapsto {\mathcal R}_D$ the
canonical projection. Let $\e_1$ be a vector field on $\Sigma_D$
such that
\begin{equation}
\label{eps1} \Pi_*\e_1(u)\in {\rm Aff}_\varphi(\lambda)\cup \bigl(-{\rm
Aff}_\varphi(\lambda)\bigr),\quad \forall u=(\lambda,\vf)\in\Sigma_D.
\end{equation}
Such fields $\e_1$
%, satisfying (\ref{eps1})
are defined modulo ${\mathcal W}_0=
%{\rm span}\,
\langle g_0,g_1,g_2\rangle$ and the sign.

The main question now is \emph{how to choose among them the canonical field, up
to the sign?}  For this first let us prove the following lemma, which will be
also useful in the sequel:

\begin{lemma}
\label{invol} The following commutative relations hold
%\begin{enumerate}
%\item $[V_{0}, {\mathcal J}]\subseteq {\mathcal J}$;
%\item $[{\mathcal J},{\mathcal J}]\subseteq {\mathcal
%J}^{(1)}$;
\begin{eqnarray}
&~&\label{i1}
 [V_i, V_i]\subseteq V_i,\quad  \forall i\geq 0;\\
&~&\label{i2} [V_i, {\mathcal J}^{(i)}]\subseteq {\mathcal
J}^{(i)}, \quad \forall i\geq 0.
\end{eqnarray}
\end{lemma}

{\bf Proof.} {\bf 1)} The proof of (\ref{i1}) is by
induction on $i$. For $i=0$ the formula is trivial, because
$V_0(\lambda)$ is the tangent space to the fiber
$(D^2)^\perp\bigl(\pi(\lambda)\bigr)$ of $(D^2)^\perp$. Now
suppose that (\ref{i1}) holds for some $i$ and prove it for
$i+1$. Take two vector fields $W_1$ and $W_2$ being tangent
to $V_{i+1}$ and prove that $[W_1, W_2]$ is tangent to
$V_{i+1}$. According to (\ref{contr2}) it is equivalent to
the fact that $\bigl[{\mathcal C},[W_1,W_2]\bigr]\subset
{\mathcal J}_{(i)}$. Note that again by (\ref{contr2}) we
have $[{\mathcal C}, W_j]\subset {\mathcal J}_{(i)}$,
$j=1,2$. Besides, by construction $V_{i+1}\subset V_i$ .
Taking into account all this, the relations (\ref{split}),
%(\ref{hV}),
and the induction hypothesis we obtain from
Jacobi identity that $$ \bigl[{\mathcal C},[W_1,W_2]\bigr]=
\bigl[[{\mathcal C},W_1],W_2]\bigr]+\bigl[W_1,[{\mathcal
C},W_2]\bigr] \subset [V_{i}\oplus {\mathcal C},
V_{i+1}]\subset {\mathcal J}_{(i)}.$$ So, $[W_1,W_2]\in
V_{i+1}$, i.e. $[V_{i+1}, V_{i+1}]\subseteq V_{i+1}$, which
completes the proof by induction of (\ref{i1}).

{\bf 2)} Directly from the definitions of ${\mathcal J}$,
$V_0$, and formula (\ref{J1}) it is easy to prove both
(\ref{i2}) for $i=0$ and the relation
\begin{equation}
\label{i3} [{\mathcal J}^{(0)},{\mathcal J}^{(0)}]\subseteq
{\mathcal J}^{(1)}.
\end{equation}

Now assume that (\ref{i2}) holds for some $i=l\geq 0$ and
prove it for $i=l+1$. Since $V_{l+1}\subset V_l$, then by
our assumption $[V_{l+1}, {\mathcal J}^{(l)}]\subset
{\mathcal J}^{(l)}$. Therefore
\begin{equation}
\label{i4} \bigl[{\mathcal C},[V_{l+1}, {\mathcal
J}^{(l)}]\bigr]\subset {\mathcal J}^{(l+1)}.
\end{equation}
On the other hand, using consequently relation (\ref{Ji1}),
the Jacobi identity, and relations (\ref{i4}), (\ref{hV}),
we obtain
\begin{multline}
\label{i5}
[V_{l+1},{\mathcal J}^{(l+1)}]=\bigl[V_{l+1}, [{\mathcal
C},{\mathcal J}^{(l)}]\bigr] \\ =\bigl[{\mathcal C},[V_{l+1},
{\mathcal J}^{(l)}]\bigr] +
\bigl[[V_{l+1},{\mathcal C}], {\mathcal J}^{(l)}\bigr]
\subset {\mathcal J}^{(l+1)}+[{\mathcal J}_{(l)},{\mathcal J}^{(l)}].
\end{multline}
If $l=0$, then ${\mathcal J}_{(0)}={\mathcal J}^{(0)}$ and
relations (\ref{i3}), (\ref{i5}) imply (\ref{i2}) for
$i=1$. If $l\geq 1$, then by (\ref{split}), the induction
hypothesis, and (\ref{Ji1})
\begin{equation}
%\begin{aligned}
\label{i6}
% ~&
[{\mathcal J}_{(l)},{\mathcal J}^{(l)}]=
[V_l\oplus {\mathcal C} ,{\mathcal J}^{(l)}]=
%\\~&
[V_l,
{\mathcal J}^{(l)}]+ [{\mathcal C},{\mathcal
J}^{(l)}]\subset {\mathcal J}^{(l+1)}.
%\end{aligned}
\end{equation}
This together with (\ref{i5}) implies that relation
(\ref{i2}) holds also for $i=l+1$. The proof of (\ref{i2})
by induction is completed. $\Box$
\medskip

 Fix some vector field $\e_1$, satisfying (\ref{eps1}).
%, satisfying (\ref{eps1}).
For any $\lambda\in \Sigma_D$ let
\begin{equation}
\label{liftJ}
{\mathcal L}_i(\lambda)=\left\{\begin{array}{ll}~\{v\in T_\lambda\Sigma_D:\Pi_*v\in\mJ_{(m-i)}(\lambda)\}&1\leq i\leq m\\
~\{v\in T_\lambda\Sigma_D:\Pi_*v\in\mJ^{(i-m)}(\lambda)\}&m\leq i\leq 2m\end{array}\right.,
\end{equation}
\begin{equation}
\label{liftV} {\mathcal W}_i(\lambda)=\{v\in
T_\lambda\Sigma_D:\Pi_*v\in V_{m-i}(\lambda)\}, 1\leq i\leq
m
\end{equation}
Then by \eqref{split}
\begin{equation}
\label{split1}
{\mathcal L}_i={\mathcal W}_i+\mathbb R h,\quad 1\leq i\leq m.
\end{equation}
According to the splitting (\ref{split1}), the projection
${\rm Pr}_{i,\lambda}$
%:{\mathcal Z}_i(\lambda)\mapsto{\mathcal W}_i(\lambda)$
from ${\mathcal L}_i(\lambda)$ onto ${\mathcal
W}_i(\lambda)$, which is parallel to $h$, is well defined.
Define the vector fields $\e_i$, $2\leq i\leq 2m$ in
addition to $\e_1$ by the following recursive formulas:

\begin{equation}
\label{basig} \e_i(\lambda)=
\left\{
\begin{array}{ll}
{\rm Pr}_{i,\lambda} \bigl([h,\e_{i-1}](\lambda)\bigr),&
2\leq i\leq m-1;\\ ~[h,\e_{i-1}], & m\leq i\leq 2m.
\end{array}
\right.
\end{equation}
Also let
\begin{equation}
\label{eta}
\eta=[\e_1,\e_{2m}].
\end{equation}

\begin{lemma}
\label{preframe}
 The tuple $(h,
\{g_i\}_{i=0}^2,\{\e_i\}_{i=1}^{2m},\eta)$ is a frame on
$\Sigma_D$.
\end{lemma}
{\bf Proof.} As before, take vector fields $E\in V_{n-4}$
and $H\in {\mathcal C}$ without stationary points and
suppose that $E$ is not collinear to the Euler field $e$.
From (\ref{basVJ}) and the fact that $\dim
\mJ^{(n-3)}=2n-4$ it follows that in order to prove the
lemma it is sufficient to prove that $[E,({\rm ad}\,
H)^{2n-7} E]\not\in \mJ^{(n-3)}$. Besides, by (\ref{i2}) we
have $[E, \mJ^{(n-4)}]\subseteq\mJ^{(n-4)}$. Therefore in
order to prove the lemma it is sufficient to prove that
$[E,\mJ^{(n-3)}]\not\subseteq \mJ^{(n-3)}$. Assuming the
converse and taking into account that
%\begin{equation}
%\label{Js}
$\mJ^{(n-3)}=\{\mathfrak s|_{(D^2)^\perp}=0\}$
%\end{equation}
 we get that
$E\in {\rm ker}\,\sigma|_{(D^2)^\perp}$, which implies that
$E$ is collinear to $H$. We have the contradiction. $\Box$

\medskip

We will say that the frame $(h,
\{g_i\}_{i=0}^2,\{\e_i\}_{i=1}^{2m},\eta)$ is {\it
associated with the vector field $\e_1$}. Note that by our
constructions
\begin{eqnarray}
&~&\label{Lspace} {\mathcal L}_i=%{\rm span}\,\bigl\{
\langle h,
\{g_i\}_{i=0}^2,\{\e_i\}_{j=1}^{i}
%\bigr\}
\rangle,\quad 0\leq i\leq
2m,\\ &~&\label{Wspace} {\mathcal W}_i=
%{\rm span}\,\bigl\{
\langle
\{g_i\}_{i=0}^2,\{\e_i\}_{j=1}^{i}
%\bigr\},
\rangle
\quad 0\leq i\leq
m-1.
\end{eqnarray}
Then by Lemma \ref{invol} we have
\begin{equation}
%\begin{alligned}
\label{kappa1} [\e_1,\e_2]=\kappa_1\e_2\,\,\,{\rm mod}\,
{\mathcal W}_1 \end{equation}
%,\quad
and, in the case $n>5$,
\begin{equation}
\label{kappa2}
[\e_1,\e_4]=\kappa_2\e_3+\kappa_3\e_4\,\,{\rm
mod}\,{\mathcal L}_2.
\end{equation}
The normalization of the field $\e_1$ can be done by
studying how the functions $\kappa_i$, $1\leq i\leq 3$, are
transformed when we pass from the frame, associated with
$\e_1$ to the frame, associated with another vector field,
satisfying (\ref{eps1}). For this first we need

\begin{lemma}
\label{gepsl} The following commutative relations hold:
\begin{eqnarray}
~&~\label{geps1}[g_1,\e_i]=(2m-2i+1)\e_i\,\,{\rm
mod}\,{\mathcal L}_{0},\\
~&~\label{geps2}
%[g_2,\e_1]\in{\mathcal L}_0,\quad
[g_2,\e_i]=(i-1)(2m-i+1)\e_{i-1}\,\, {\rm mod}\,{\mathcal
L}_{0},
%\,\,i\geq 2
\\~&~\label{geps3}[g_0,\e_i]=-\e_i
\,\,{\rm mod}\,{\mathcal L}_0 \end{eqnarray} In the case
$1\leq i\leq m-1$ the subspace ${\mathcal L}_0={\rm
  span}\{h,g_0,g_1,g_2\}$ can be
replaced by ${\mathcal W}_0={\rm span}\{g_0,g_1,g_2\}$.
\end{lemma}
{\bf Proof.} Let $\psi:\mathbb R\mapsto \mathbb R$,
$\psi_0=0$. Then from (\ref{norme1}) it follows easily that
\begin{equation}
\label{repare1}
\e_{\psi\circ\vf}(\lambda)=(\psi'(0))^{m-\frac{1}{2}}\e_\vf(\lambda)\,\,{\rm
mod} \{{\mathbb R} e(\lambda)\}.
\end{equation}
Using the last identity and (\ref{flows}) one gets without
difficulties that \begin{eqnarray} ~&~\label{Ad1}
(F_{1,s}^{-1})_*(\e_1)=e^{(2m-1)s}\e_1 \,\,\,\,{\rm
mod}\,\, {\mathcal W}_0,
%{\rm span}\{g_0,g_1,g_2\}
\\ ~&~ \label{Ad2}\
(F_{2,s}^{-1})_*(\e_1)=\e_1 \,\,\,\,{\rm mod}\,\, {\mathcal
W}_0.
%{\rm span}\{g_0,g_1,g_2\}
\end{eqnarray}
In addition, let us show that   %Using in addition (\ref{homoth}) and
%(\ref{flow0}), one gets
\begin{equation}
\label{Ad0} (F_{0,s}^{-1})_*(\e_1)=e^{-s}\e_1 \,\,\,\,{\rm
mod}\,\,
{\mathcal W}_0
%{\rm span}\{g_0,g_1,g_2\}
\end{equation}
Indeed, the homotheties $\delta_s$, defined by
(\ref{homoth}), preserve both the characteristic line
distribution ${\mathcal C}$ and the distribution ${\mathcal
J}$. Therefore it preserves also all ${\mathcal J}^{(i)}$.
Note also that
\begin{equation}
\label{deltas} \delta_s^*\sigma=e^s\sigma, \end{equation}
which implies that $\delta_s$ preserve all ${\mathcal
J}_{(i)}$ as well. This together with (\ref{flow0}) yields
that there exists a function $\alpha$ such that
$(F_{0,s}^{-1})_*(\e_1)=\alpha \e_1 \,\,\,\,{\rm mod}\,\,
{\mathcal W}_0$. The last formula together with
(\ref{basig}) implies that
\begin{equation}
\label{F0ei} (F_{0,s}^{-1})_*(\e_i)=\alpha \e_i
\,\,\,\,{\rm mod}\,\, {\mathcal L}_{i-1},\quad 1\leq i\leq
2m.
\end{equation}
From the normalization (\ref{norme1}) of the vector field
$\e_1$ it follows that
\begin{equation}
\label{norme11} \Pi^*\sigma(\e_{m+1},\e_m)=1,
\end{equation}
which together with (\ref{F0ei}) implies that
\begin{equation}
\label{PF1} (\Pi\circ
F_{0,s}^{-1})^*\sigma(\e_{m+1},\e_m)=\alpha^2.
\end{equation}
On the other hand, by (\ref{flow0}), $\Pi\circ
F_{0,s}^{-1}=\delta_{2s}^{-1}\circ\Pi$, which together with
(\ref{deltas}) implies that
\begin{equation}
\label{PF2} (\Pi\circ
F_{0,s}^{-1})^*\sigma(\e_{m+1},\e_m)=e^{-2s}.
\end{equation}
Combining (\ref{PF1}) and (\ref{PF2}), we get
$\alpha=e^{-s}$, which proves (\ref{Ad0}).

Relations (\ref{Ad1}),(\ref{Ad2}), and (\ref{Ad0}) imply
(\ref{geps1}), (\ref{geps2}), and (\ref{geps3}) respectively for
$i=1$. Then relations (\ref{geps1})-(\ref{geps3}) for $i>1$ can be
proved by induction, using (\ref{gl2}), (\ref{basig}), and the Jacobi
identity. The last sentence of the lemma follows from the fact that by
our constructions the vectors $\pi_*\circ\Pi_*(\e_i)=0$ for $1\leq
i\leq m-1$.$\Box$
\medskip

The following lemma gives the normalization of the fields
$\e_1$:

\begin{lemma}
\label{normel} Among all vector fields , satisfying
(\ref{eps1}), there exists a unique, up to the sign, field
$\tilde\e_1$ such that all functions $\kappa_i$, $1\leq
i\leq 3$,
%, appearing in (\ref{kappa}),
are identically zero, namely, the following commutative
relations hold \begin{equation} \label{normrel} [\tilde
\e_1, \tilde \e_2]\in{\mathcal W}_1,\quad [\tilde
\e_1,\tilde \e_4]\in{\mathcal L}_2.
\end{equation}

\end{lemma}
{\bf Proof.} Let $\e_1$ and $\tilde\e_1$ be two vector
field, satisfying (\ref{eps1}). Then there exist functions
$\{\mu_i\}_{i=0}^2$ such that
\begin{equation*}
%\label{mu}
\tilde \e_1 =\pm\e_1+\mu_0 g_0+\mu_1 g_1+\mu_2
g_2.
\end{equation*}
Using (\ref{gl2}) it is easy to show that
$[h,\tilde\e_2]=\pm[h,\e_2] +2\mu_1 h\,\,{\rm mod}\,
{\mathcal W}_0$, which implies that $\tilde e_2=e_2 {\rm
mod}\,{\mathcal W}_0$. In the same way one can show that
\begin{equation}
\label{ee} \begin{aligned} ~&\e_i=\pm\e_i\,\,{\rm
mod}\,{\mathcal W}_0,\quad 2\leq i\leq m-1;\\ ~&\tilde
\e_i=\pm\e_i\,\,{\rm mod}\,{\mathcal L}_0,\quad m\leq i \leq 2m.
\end{aligned}
\end{equation}
Suppose that the frame, associated with $\e_1$, satisfies
relations (\ref{kappa1}) and (\ref{kappa2}), while the
frame, associated with $\tilde\e_1$ satisfies
\begin{eqnarray*}
%\begin{alligned}
~&~
[\tilde\e_1,\tilde\e_2]=\tilde\kappa_1\tilde\e_2\,\,\,{\rm
mod}\, {\mathcal W}_1\nonumber\\
%\end{equation}
%,\quad
%and, in the case $n>5$,
%\begin{equation}
~&~
[\tilde\e_1,\tilde\e_4]=\tilde\kappa_2\tilde\e_3+\tilde\kappa_3\tilde\e_4\,\,{\rm
mod}\,{\mathcal L}_2\nonumber .
%\end{equation}
\end{eqnarray*}
Then by direct computation, using relations
(\ref{geps1})-(\ref{geps3}) and (\ref{ee}), one can show
without difficulties that \begin{equation} \label{transk}
\left\{\begin{aligned} ~&
\tilde\kappa_1=\kappa_1\pm(2m-3)\mu_1-\mu_0\\
~&\tilde\kappa_2=\kappa_2\pm(6m-9)\mu_2\\
~&\tilde\kappa_3=\kappa_3\pm(2m-7)\mu_1-\mu_0
\end{aligned} \right. .
\end{equation}
Obviously, system (\ref{transk}) w.r.t. $\mu_i$, $0\leq
i\leq 2$, has a unique, up to the sign, solution, when
$\tilde\kappa_i=0$, $1\leq i\leq 2$, which completes the
proof of the lemma. $\Box$
\medskip

Suppose that $\tilde\varepsilon_1$ is one of the two vector
fields, found in the previous lemma. Then two frames $(h,
\{g_i\}_{i=0}^2,\{\tilde\e_i\}_{i=1}^{2m},\eta)$ and $(h,
\{g_i\}_{i=0}^2,\{-\tilde\e_i\}_{i=1}^{2m},\eta)$ are
canonically defined on $\Sigma_D$. They are {\it the
canonical frames} we were looking for. Also, this immediately
implies that the groups of symmetries of
$(2,n)$-distributions, $n>5$, of maximal class is at most
$(2n-1)$-dimensional. The proof of the theorem is
completed. $\Box$

\begin{remark}
{\rm The normalization, implemented above for $n>5$, does
not work in the case $n=5$. In this case $m=2$ and the
relation (\ref{kappa2}) is not true (actually, in this case
$[\e_1,\e_4]=\eta$). Note also that the bracket
$[\e_1,\e_3]$ cannot give new conditions with compare to
the bracket $[\e_1,\e_2]$,
because
\begin{equation}
\label{e1e3} [\e_1,\e_3]=\bigl[h,[\e_1,\e_2]\bigr].
\end{equation}
 It is
also in accordance with the result of E.~Cartan~\cite{cartan}.
In the case $n=5$ the canonical frame does not exist on
$\Sigma_D$ (which is $9$-dimensional), one has to prolong
further to construct it.} $\Box$
\end{remark}

\begin{remark}
\label{newnorm} {\rm One can suggest normalizations of the
vector field $\e_1$ different from one given by Lemma
\ref{normel}. For example, among all vector fields ,
satisfying (\ref{eps1}) there exist unique, up to the sign,
vector field $\e_1$ such that $[h,\e_1]\in {\mathcal W}_2$,
$[\e_1,\e_2]\in {\mathcal W}_1 $, and $[\e_1,\e_4]\in{\rm
span}\,\{\e_4,{\mathcal L}_2\}$. In the case $n>6$ we can
get one more normalization by replacing the last condition
of the previous normalization by $[h,\e_2]\in {\mathcal
W}_3$. The frames, associated with such $\e_1$, are also
intrinsically defined, up to the corresponding reflection.
In particular, Theorem \ref{equivprob} below remains true,
if one uses these frames instead of the canonical frames,
constructed in the proof of Theorem \ref{frtheor}. Note
also that both of these normalizations cannot be
implemented in the case $n=5$, because in this case the
vector $[h,\e_1]$ is not tangent to the fiber of
$\Sigma_D$, considered as the fiber bundle over $M$.$\Box$
}\end{remark}

\begin{theor}
\label{equivprob} Two rank 2 distributions $D_1$ and $D_2$ of maximal class are
equivalent iff there exists the diffeomorphism ${\mathcal
F}:\Sigma_{D_1}\mapsto\Sigma_{D_2}$, which transform one of the canonical frames of $D_1$ to one of the canonical
frames of $D_2$.
%%fields $\vec e$, $\overrightarrow  {h_A}$, $\varepsilon_A$ to the
%%following three vector fields $\vec e$, $\overrightarrow {\bar
%%h_A}$, $\bar\varepsilon_A$
%canonical frame of $\bar D$,
%%i.e.,
%%\begin{equation}
%%\label{fbarf}
%,\,\,\,\, \forall \lambda\in \aleph_D,\,\,\,
%1\leq i\leq 7.
%% \begin{array} {c}
%% {\mathcal F}_*\vec e(\lambda)= \vec e\Bigr({\mathcal
%%F}(\lambda)\Bigr) \\ {\mathcal F}_*\varepsilon_A(\lambda) =
%%\bar\varepsilon_A \Bigr({\mathcal F}(\lambda)\Bigr)\\ {\mathcal F}_*
%%\overrightarrow  {h_A}(\lambda)= \overrightarrow {\bar
%%h_A}\Bigr({\mathcal F}(\lambda)\Bigr)
%%\end{array}
%%\end{equation}
\end{theor}

{\bf Proof.} The necessity is obvious.
%ecause the construction of the canonical frame is intrinsic.
%f the canonical frame of the distribution.
Let us prove the
sufficiency.
Let $(h^k,
\{g_i^k\}_{i=0}^2,\{\tilde\e_i^k\}_{i=1}^{2m},\eta^k)$ are the canonical frames of the distributions $D_k$
respectively, where $k=1,2$, such that ${\mathcal F}$ transforms the frame with $k=1$ to the frame with $k=2$.
If the maps $\Pi_k:\Sigma_{D_k}\mapsto {\mathcal R}_{D_k}$
 and $\pi:T^*M\mapsto M$ are the canonical projections, then the map
 ${\mathfrak p}_k\stackrel{def}{=}\Pi_k\circ\pi$ defines the fiber bundle $\Sigma_{D_k}$ over $M$.
 By our constructions,  the tangent spaces to this fibers coincide with
 $%{\mathcal W}_{m-1}=
 {\rm span}\{g_0,g_1,g_2,\e_1^k,\ldots, \e_{m-1}^k\}$.
 %The pairs of vector fields $\vec e$, $\varepsilon_A$
%and $\vec e$, $\bar\varepsilon_A$ constitute the bases of the
%fibers of $(D^2)^\perp$ and $(\bar D^2)^\perp$ respectively.
This
%together with the first two relations of (\ref{fbarf})
yields
that the diffeomorphism ${\mathcal F}$ is fiberwise, i.e., there
exists the diffeomorphism $F:M\mapsto M$ such that
\begin{equation}
\label{fiber}
 F\circ {\mathfrak p}_1={\mathfrak p}_2\circ {\mathcal F}.
\end{equation}
Note that by our constructions $\Pi_* h^1$ and $\Pi_* h^2$
span the characteristic line distributions of $D$ and $\bar
D$ respectively. Therefore from Remark \ref{projCr} and
relations (\ref{fiber}), (\ref{projC}) it follows easily
that
% Note also that for any
%$q\in M$ $${\rm span}\Bigr(\pi_*\Bigl(\{
%\overrightarrow{h_A}(\lambda):\lambda\in\aleph_D,\pi(\lambda) =q\}
%\Bigr)\Bigr)=D(q),$$
%%(recall that $\wp_2= \overrightarrow {h_A}$)
%$${\rm span}\Bigl(\pi_*\Bigl(\{ \overrightarrow{\bar h_A}
%(\lambda):\lambda\in\aleph_{\bar D},\pi(\lambda)=q\}
%\Bigr)\Bigr)=\bar D(q).$$ Using the last two relations, the third
%relation of (\ref{fbarf}), and (\ref{fiber})
$F_* D=\bar D$, which means that the distributions $D$ and
$\bar D$ are equivalent. $\Box$
\medskip

Now we will list several properties of the canonical
frames.
%, introduced in the proof of the Theorem 1.
First from (\ref{normrel}) and (\ref{e1e3}) it
follows that
\begin{equation}
\label{te1e3} [\tilde\e_1,\tilde\e_3]\in{\mathcal L}_2.
\end{equation}
Second, directly from Lemma \ref{invol} it follows that
\begin{eqnarray}
&~&\label{bas1}
[\tilde\varepsilon_1,\tilde\varepsilon_i]\in {\mathcal
L}_i,\quad 5\leq i\leq 2n-6;\\ \label{bas2}
&~&[\tilde\varepsilon_{i_1},\tilde\varepsilon_{i_2}]\in{\mathcal
L}_{i_2+1},\quad 2\leq i_1\leq i_2\leq 2n-6-i_1
\end{eqnarray}
where the subspaces ${\mathcal L}_i$ are as in (\ref{Lspace}).
Further, from identities (\ref{last0}), (\ref{A2m-2}), and
(\ref{A2m-3})
%of Remark \ref{projcrit}
it follows that
\begin{equation}
\label{hlast} [h, \tilde\varepsilon_{2n-6}]\in {\mathcal
L}_{2n-9}. \end{equation} Finally we have
\begin{lemma}
The following relation holds
\begin{equation}
\label{heisen}
 [\tilde\varepsilon_i,\tilde\varepsilon_{2n-5-i}]=(-1)^{i+1}\eta\quad
 {\rm mod}\,\, {\mathcal
L}_{2n-6}, \quad 1\leq i\leq n-3.\end{equation} \end{lemma}

{\bf Proof.} The proof is by induction on $i$. For $i=1$
formula (\ref{heisen}) follows from (\ref{eta}). Now assume
that (\ref{heisen}) holds for some $i=j$ and prove it for
$i=j+1$. By (\ref{bas2})
$$[\tilde\varepsilon_j,\tilde\varepsilon_{2n-6-j}]\in
{\mathcal L}_{2n-5-j}$$ Taking Lie brackets with $h$ from
both sides of the last inclusion and using formulas
(\ref{basig}) and (\ref{hlast}) together with the Jacobi
identity, one gets $$[\tilde\varepsilon_{j+1},
\tilde\varepsilon_{2n-5-(j+1)}]+[\tilde\varepsilon_j,\tilde\varepsilon_{2n-6-j}]\in
{\mathcal L}_{2n-6}.$$ Therefore by induction hypothesis
for $i=j$ $$ [\tilde\varepsilon_{j+1},
\tilde\varepsilon_{2n-5-(j+1)}]=-[\tilde\varepsilon_j,\tilde\varepsilon_{2n-6-j}]\,\,
{\rm mod}\,\, {\mathcal L}_{2n-6}=(-1)^{j+2} \eta \,\, {\rm
mod}\,\, {\mathcal L}_{2n-6},$$ which proves formula
(\ref{heisen}) also for $i=j+1$. The proof by induction is
completed. $\Box$

\section{The most symmetric case}
\setcounter{equation}{0} \indent

The present section is devoted to the following

\begin{theor}
\label{symtheor} Let $n>5$. Then any $(2,n)$-distribution of maximal
class with $(2n-1)$-dimensional Lie algebra of infinitesimal
symmetries is locally equivalent to the distribution, associated
with the underdetermined ODE $z'(x)=\bigl(y^{(n-3)}(x)\bigr)^2$. The
symmetry algebra of this distribution is isomorphic to a semidirect
sum of $\mathfrak{gl}(2,\mathbb R)$ and $(2n-5)$-dimensional
Heisenberg algebra ${\mathfrak n}_{2n-5}$.
\end{theor}

{\bf Proof.}
%First, we will prove that any $(2,n)$-distribution
%of maximal class with $(2n-1)$-dimensional group of
%symmetries (if it exists) has algebra of infinitesimal
%symmetries, which is isomorphic to a semidirect sum of
%$\mathfrak{gl}(2,\mathbb R)$ and $(2n-5)$-dimensional
%Heisenberg algebra ${\mathfrak n}_{2n-5}$.
If a $(2,n)$-distribution of maximal class has a
$(2n-1)$-dimensional group of symmetries, then all
structural functions of its canonical frames have to be
constant. Using this fact we get the following

\begin{lemma}
\label{all0}
In addition to (\ref{gl2}),
the only nonzero commutative
relations of each of the canonical frames of a $(2,n)$-
distribution with a
$(2n-1)$-dimensional group of symmetries
 are
\begin{equation}
\label{commrel}
\begin{aligned}
~&[h,\tilde\e_i]=\tilde\e_{i+1},\,\,[\tilde \e_i,\tilde
\e_{2m-i+1}]=(-1)^{i+1}\eta,\,\, [g_1, \tilde
\e_i]=(2m-2i+1)\tilde\e_i,\\ ~&[g_2,\tilde
\e_i]=(i-1)(2m+1-i)\tilde\e_{i-1},\,\,[g_0,
\tilde\e_i]=-\tilde\e_i, \,\, [g_0,\eta]=-2\eta.
\end{aligned}
\end{equation}
\end{lemma}

{\bf Proof} By Lemma \ref{gepsl}
\begin{equation}
\label{const1}
[g_0,\tilde\e_1]=-\tilde\e_1+\sum_{i=0}^2\alpha_i g_i.
\end{equation}
 Let us
prove that $\alpha_i=0$ for all $0\leq i\leq 2$. Indeed,
from (\ref{normrel}) and Lemma \ref{gepsl} it follows that
\begin{equation}
\label{const2}
\bigl[g_0,[\tilde\e_1,\tilde\e_2]\bigr]\in{\mathcal
W}_1,\quad\bigl[g_0,[\tilde\e_1,\tilde\e_4]\bigr]\in{\mathcal
L}_2
\end{equation}
%By (\ref{basig}), (\ref{gl2}), and Jacobi identity one gets
%easily that $$[g_0,\tilde\e_2]=-\tilde e_2+2\alpha_1
%h-\alpha_2 g_1,\quad [g_0,\tilde\e_4]=-\tilde\e_4.$$
Using the Jacobi identity and relation \eqref{geps3} we get easily
that
\begin{equation*}
\begin{aligned}
~&[g_0,[\tilde\e_1,\tilde\e_2]=\bigl((2m-3)\alpha_1-\alpha_0\bigr)\tilde\e_2\,\,
{\rm mod}\,{\mathcal W}_1,\\ ~& [g_0,[\tilde\e_1,\tilde\e_4]]=
\bigl((2m-7)\alpha_1-\alpha_0\bigr)\tilde\e_4+(6m-9)\alpha_2\tilde
\e_3\,\, {\rm mod}\,{\mathcal L}_2. \end{aligned}
\end{equation*} Comparing the last relations with
(\ref{const2}), we have immediately that $\alpha_i=0$ for
all $0\leq i\leq 2$. In other words,
$[g_0,\tilde\e_1]=-\tilde\e_1$. Let us prove that
%This together with (\ref{basig}) implies in turn that
\begin{equation}
\label{const3} [g_0,\tilde\e_i]=-\tilde\e_i,\quad 1\leq
i\leq 2m.
\end{equation}
The proof is by induction. For $i=1$ relation
(\ref{const3}) is true. Suppose that it is true for some
$i=j$ and prove it for $i=j+1$. By (\ref{basig})
\begin{equation}
\label{const31} [h,\tilde\e_j]=\tilde\e_{j+1}+\tau_j h,
\end{equation}
where $\tau_j\equiv 0$ for $m-1\leq j\leq 2m-1$. Taking the
Lie brackets with $g_0$ from both sides of the last
identity and using (\ref{gl2}), the induction hypothesis,
and (\ref{const31}) again, one gets easily that
\begin{equation}
\label{const32}
[g_0,\tilde\e_{j+1}]=-[h,\e_j]=-\tilde\e_{j+1}-\tau_j h.
\end{equation}
If $m-1\leq j\leq 2m-1$ we get (\ref{const3}) for $i=j+1$,
because $\tau_j=0$. If $1 \leq j\leq m-2$, then from Lemma
\ref{gepsl} it follows that
$[g_0,\tilde\e_{j+1}]\in{\mathcal W}_{j+1}$. Therefore
(\ref{const32}) implies that $\tau_j=0$. Hence
(\ref{const3}) holds again. The proof by induction of
(\ref{const3}) is completed. Note that we have proved at
the same time that
\begin{equation}
\label{basig1} [h,\tilde\e_i]=\tilde\e_{i+1},\quad 1\leq
i\leq 2m-1
\end{equation}
Further, from (\ref{const3}) and the Jacobi identity it
follows immediately that
\begin{equation}
\label{const35}
[g_0,\eta]=\bigl[g_0,[\tilde\e_1,\tilde\e_{2m}]\bigr]=-2\eta.
\end{equation}
The identity (\ref{const3}) allows to show that a lot of structural
constants of the canonical frame vanish. First, by Lemma \ref{gepsl}
$$[g_1,\tilde \e_1]=(2m-1)\tilde \e_1+\sum_{i=0}^2\beta_i g_i.$$
Taking Lie brackets with $g_0$ from both sides of the last identity
and comparing the coefficients with the help of (\ref{const3}),
(\ref{gl2}), and the Jacobi identity, one obtains immediately that
$\beta_i=0$ for all $0\leq i\leq 2$, which together with
(\ref{basig1}) implies in turn that
\begin{equation}
\label{const4} [g_1,\tilde\e_i]=(2m-2i+1)\tilde\e_i,\quad 1\leq
i\leq 2n-6.
\end{equation}
The last identity yields also that
\begin{equation}
\label{const5}
[g_1,\eta]=\bigl[g_1,[\tilde\e_1,\tilde\e_{2m}]\bigr]=0.
\end{equation}

In the same way we get
\begin{equation}
\label{const6} [g_2,\tilde\e_i]=(i-1)(2m-i+1)\e_i,\quad
1\leq i\leq 2n-6.
\end{equation}
Further, suppose that $$[\tilde\e_i,\tilde
\e_j]=\sum_{k=1}^{2n-6}a_{ij}^k\tilde\e_k+\sum_{k=0}^2
b_{ij} ^k g_k+c_{ij} h+d_{ij}\eta.$$ Again taking Lie
brackets with the field $g_0$ from both sides and comparing
the coefficients with the help of (\ref{geps3}),
(\ref{gl2}) and the Jacobi identity we get immediately
that $a_{ij}^k=0$, $b_{ij}^k=0$ and $c_{ij}^k=0$. In other
words, $[\tilde\e_i,\tilde \e_j]=d_{ij}\eta$. Taking Lie
brackets with $g_1$ from both sides of the last identity
and comparing the coefficients with the help of
(\ref{const4}) and the Jacobi identity one obtains easily
that $d_{ij}=0$ for $i+j\neq 2m+1$. On the other hand, by
(\ref{heisen}) we have $d_{ij}=(-1)^{i+1}$ for $i+j=2m+1$.
In other words,
\begin{equation}
\label{const7}
[\tilde\e_i,\tilde\e_j]=\left\{\begin{array}{ll} (-1)^{i+1}\eta & i+j=2m+1\\
0&i+j\neq 2m+1
\end{array}\right..
\end{equation}
This together with the Jacobi identity immediately implies
that
$[\tilde\e_i,\eta]=\bigl[\tilde\e_i,[\tilde\e_1,\tilde\e_{2m}]\bigr]=0$
for $1<i<2m$. To prove that $[\tilde \e_i,\eta]=0$ also for
$i=1$ and $i=2m$ we take Lie brackets with $g_0$ and
compare the coefficients with the help of
(\ref{const3}),(\ref{const35}), and the Jacobi identity.
Finally, to prove that $[h,\tilde \e_{2m}]=0$ we take Lie
brackets with $g_1$ and compare the coefficients with the
help of (\ref{const4}) and the Jacobi identity. This
completes the prove of the lemma. $\Box$
\medskip

The previous lemma and Theorem \ref{equivprob} imply the
uniqueness, up to the equivalence, of the germ of
$(2,n)$-distribution of maximal class with
$(2n-1)$-dimensional group of symmetries. Besides, from
these relations it follows that the algebra of
infinitesimal symmetries of such distribution is isomorphic
to the semi-direct sum of ${\mathfrak{gl}}(2,{\mathbb R})$
($\sim{\rm span}_{\mathbb R}\,\{ g_0,g_1,g_2,h\}$) and the
Heisenberg group ${\mathfrak n}_{2m+1}$ ($\sim {\rm
span}_{\mathbb R}\{\tilde
\e_1,\ldots,\tilde\e_{2m},\eta\}$).
%Besides, from the first
%of these relations it
%follows that %the linear span (over ${\mathbb R}$) of the
%%vector fields
%${\rm span}){\mathbb R}(\tilde
%\e_1,\ldots,\tilde\e_{2m},\eta)$ is
%%a Lie algebra
%isomorphic to the Heisenberg group ${\mathfrak n}_{2m+1}$,
%while all other relations together with the relations
%$[h,\tilde \e_i]=\tilde \e_{i+1}$, $1\leq i\leq 2m$ define
%the action of ${\rm span}_{\mathbb R}\,\{ g_0,g_1,g_2,h\}$
%($\sim{\mathfrak{gl}}(2,{\mathbb R})$) on ${\mathfrak
%n}_{2m+1}$.
%%this Heisenberg
%%group.

To complete the proof of the theorem it remains to show
%by direct computations one can show
%easily
that for the $(2,n)$-distribution $D_0$, associated with
the underdetermined ODE
$z'(x)=\frac{1}{2}\bigl(y^{(n-3)}(x)\bigr)^2$, the only
nonzero commutative relations for its canonical frames are
(\ref{gl2}) and (\ref{commrel}). Let $p_i=y^{(i)}$, as in
Introduction. Then the distribution $D_0$ is given in
$M={\mathbb R}^n$ with coordinates $(x,p_0,\ldots,
p_{n-3},z)$ by the intersection of the annihilators of the
forms
\begin{equation}
\label{Pfaff1}
\begin{aligned}
~& dp_i-p_{i+1} dx , \,\,0\leq i\leq n-4,\\
%\quad
% dq_j-q_{j+1} dx,\,\,
%0\leq j\leq 0\leq n-4,\\
~& dz-\frac{1}{2}p_{n-3}^2\,dx
 %F(x,p_0,\ldots, p_s,q_0,\ldots, q_{r-1})dx.
\end{aligned}
\end{equation}
Set, as before, $m=n-3$. From (\ref{Pfaff1}) the following
two vector fields span the distribution $D_0$:
\begin{equation} \label{basf} X_1=\partial_{p_m},\quad X_2=\partial_x+
\sum_{i=0}^{m-1}
p_{i+1}\partial_{p_i}+\frac{1}{2}p_m^2\partial_z.
\end{equation}
Let vector fields $X_3$, $X_4$, $X_5$ are as in
(\ref{x345}). Then
\begin{equation}
\label{x345f}
X_3=\partial_{p_{m-1}}+p_m\partial_z,\,\,X_4=\partial_z,\,\,
X_5=-\partial_{p_{m-2}}.
\end{equation}
Denote
\begin{equation}
\label{xn} X_i=({\rm ad} X_2)^{i-5}X_5=(-1)^i\partial_
{p_{m+3-i}},\quad 6\leq i\leq m+3 \end{equation} Let
$u_i:T*M\mapsto{\mathbb R}$, $1\leq i\leq m+3$, be the corresponding
quasi-impulses, defined by (\ref{quasi25}). Then the tuple
$(x,p_0,\ldots, p_m,z,u_1,\ldots, u_{m+3})$ defines the coordinates
on $T^*M$. It is clear that in this coordinates the submanifold
$(D^2)^\perp$ is given by
equations $(D_0^2)^\perp=\{u_1=u_2=u_3=0\}$.% Therefore the
%tuple $(x,p_0,\ldots, p_m,z,u_4,\ldots, u_{m+3}) $defines
%the coordinates on $(D_0)^2$. In the sequel we will work

 Denote by $\overline X_i$ the vector field on
$(D_0^2)^\perp$, which is the lift of the vector field $X_i$ (i.e.,
$\pi_* \overline X_i=X_i$) such that $du_j(\overline X_i)=0$ for all
$1\leq j\leq n$. In addition to (\ref{x345}) the only nonzero
commutators generated by the vector fields $\{X_i\}_{i=1}^5$ are
$({\rm ad} X_2)^j X_5$, $1\leq m-2$. This together with
\eqref{foli25} implies that the characteristic line distribution
${\mathcal C}$ satisfies $${\mathcal C}= \bigl\{\mathbb{R} \bigl(u_4
\overline X _2-u_5\overline X_1+u_4\sum_{i=5}^{m+2}
u_{i+1}\partial_{u_i}\bigr)\bigr\}.$$ From this it is not difficult
to show that
\begin{equation}
\label{RD0} {\mathcal R}_{D_0}=\{\lambda\in (D_0^2)^\perp:
u_4(\lambda)\neq 0\}
\end{equation}
(if $u_4(\lambda)=0$, then
$\mJ^{(3)}(\lambda)=\mJ^{(2)}(\lambda)$). Define the
following vector field $H$ on ${\mathcal R}_{D_0}$, which
generates the characteristic line distribution ${\mathcal
C}$:
\begin{equation} \label{homH} H= \overline X
_2-\frac{u_5}{u_4}\overline X_1+\sum_{i=5}^{m+2}
u_{i+1}\partial_{u_i} \end{equation} Assume that $\vf_\lambda$ is
the parameterization of the characteristic curve $\gamma$, passing
through $\lambda$, such that
\begin{equation}
\label{tt} \vf_\lambda(e^{tH}\lambda)=t.
\end{equation}
Then by direct computations one can show that the following vector
$$\e_{\vf_\lambda}(\lambda)=|u_4(\lambda)|^{1/2}\partial_{u_{m+3}}(\lambda)$$
%from Lemma \ref{praglem}
satisfies
\begin{equation} \label{evf0}
 \e_{\vf_\lambda}(\lambda)\in {\rm Aff}_{\vf_\lambda}(\lambda)\cup \bigl(-{\rm
Aff}_{\varphi_\lambda}(\lambda)\bigr)
 %=|u_4(\lambda)|^{1/2}\partial_{u_{m+3}}(\lambda).
 \end{equation}
 Denote by $\e_H$ the vector field, satisfying $\e_H (\lambda)=
 \e_{\vf_\lambda}(\lambda)$ for all $\lambda\in {\mathcal
 R}_{D_0}$.
 By direct calculation it is easy to show that
 \begin{equation}
 \label{ad}
 \begin{aligned}
 ~&({\rm ad} H)^i\e_H=
 (-1)^i|u_4|^{1/2}\partial_{u_{m+3-i}},\quad 0\leq i\leq
 m-2,\\
 ~&({\rm ad} H)^{m-1}\e_H=(-1)^{m-2}|u_4|^{-1/2}\overline
 X_1,\quad ({\rm ad} H)^m\e_H =(-1)^{m-1}|u_4|^{-1/2}\overline
 X_3,\\
 ~&({\rm ad} H)^{m+i}\e_H=(-1)^{m-1}|u_4|^{-1/2}\bigl(\overline
 X_{m+i}-\frac{u_{4+i}}{u_4}\overline X_4\bigr),\quad 1\leq
 i\leq  m-1,
 \end{aligned}
 \end{equation}
 and finally
 \begin{equation}
 \label{flatH}
 ({\rm ad} H)^{2m}\e_H=0.
 \end{equation}
 Then by \eqref{A2m-2}
 %Remark \ref{projcrit}
 the parameterizations
 $\vf_\lambda$ defined by (\ref{tt}) are projective\footnote{
 Actually from (\ref{flatH}) it follows that the curves in a projective space
 associated with all abnormal extremals  of the distribution $D_0$
 %parametrized with the help of
 %the field $H$, are so-called flat curves in the corresponding
 %Lagrangian  Grassmannians (see \cite{jac2} for the definition
 %of the flat curves).
are so-called normal rational curves.
 The interesting question is whether the distribution $D_0$ is
 a unique, up to the equivalence, $(2,n)$-distribution, having this
 property.}; namely,
 $\vf_\lambda(\cdot)\in{\mathfrak P}_\lambda$. Therefore if
 some parameterization $\bar\vf$ belongs to ${\mathfrak
 P}_\lambda$, then
 $\bar\vf=\frac{a\vf_\lambda}{b\vf_\lambda+1}$. Let us
 introduce the coordinates on $\Sigma_{D_0}$ in the
 following way:
 $$\left(\lambda,\frac{a\vf_\lambda}{b\vf_\lambda+1}\right)\mapsto
 (x, p_0,\ldots , p_m, x,z,u_4,\ldots, u_n, a, b)$$
 Then  from (\ref{flows}), (\ref{flow0}), and (\ref{flowh}
 it is not difficult to show that
 \begin{equation}
 \label{gl2coord}
 \begin{aligned}
 ~&g_1=2a\partial_a,\quad g_2=-a\partial_b,\quad g_0=2\sum_{i=4}^n
 u_i\partial_{u_i}\\
 ~&h=\overline H-2b\partial_a-\frac{b^2}{a}\partial_b,
 \end{aligned}
 \end{equation}
 where   the vector field $\overline H$ is the lift of the vector field $H$
 on $\Sigma_{D_0}$ (i.e., $\Pi_*
\overline H=H$) such that $da(H)=db(H)=0$. From
(\ref{repare1}) and (\ref{evf0}) it follows that $\e_{\bar
vf}=a^{m-1/2} \e_H\,\,{\rm mod} \{{\mathbb R} e\}$, where
$\bar\vf=\frac{a\vf_\lambda}{b\vf_\lambda+1}$ and $e$ is
the Euler field on $T^*M$. Therefore the vector field
\begin{equation}
\label{e1flat} \e_1=a^{m-1/2}|u_4|^{1/2}\partial_{u_{m+3}}
 \end{equation}
satisfies (\ref{eps1}). Consider the frame $(h,
\{g_i\}_{i=0}^2,\{\e_i\}_{i=1}^{2m},\eta)$ associated with
the vector field $\e_1$. Using commutative relations
(\ref{ad}) and formulas (\ref{gl2coord}), (\ref{e1flat}),
it is easy to show that this frame is canonical and the
only its nonzero commutative relations are (\ref{gl2}) and
(\ref{commrel}). This completes the proof of the theorem.
 $\Box$
%\begin{remark}
%\label{open} {\rm
% From Remark 3.4 of \cite{zelvar} it
%follows that a rank 2 distribution has minimal class $1$ at
%a point $q$ if and only if $\dim\, D^3(q)=4$. On the other
%hand if a distribution has small growth vector
%$(2,3,4,\ldots)$ on some open set, then one can get from
%this distribution a distribution with small growth vector
%$(2,3,5, \ldots)$. In other words The natural question is
\medskip

\begin{remark}
{\rm Note that for the most symmetric case the frames,
indicated in Remark \ref{newnorm} coincide (up to the
reflection) with the canonical frames, introduced in the
proof of Theorem \ref{frtheor}.} $\Box$
\end{remark}

\section{Discussion}
\subsection{Distributions of non-maximal class}
As was already mentioned before, a rank~2 distribution $D$
has the smallest possible class $1$ at a point $q$ iff
$\dim\, D^3(q)=4$ (see \cite[Remark~3.4]{zelvar}). Suppose
that $D$ satisfies $\dim\, D^3(q)=4$ on some open set
$M^o$. It is easy to show that the distribution $D^2$ has a
one-dimensional characteristic distribution $C$. Then
(locally) we can consider the quotient $M'$ of the manifold
$M^o$ by the corresponding one-dimensional foliation
together with a new rank~2 distribution $D'$ obtained by
the factorization of $D^2$.

In fact, $D$ can be uniquely reconstructed from $D'$. Let $P(D')$ be a
submanifold in $P(TM')$ consisting of all lines lying in $D'$.
Similarly to the canonical contact system on $P(TM')$, we can define
lifts of integral curves of $D'$ to $P(D')$ and a canonical rank 2
distribution on $P(D')$ generated by tangent vectors to these lifts.
It can be proved that this contact system on $P(D')$ is locally
equivalent to $D$.

Iterating this procedure, we end up either at a nonholonomic rank 2
distribution on a three-dimensional manifold or at a distribution
$\widetilde D$, satisfying $\dim \widetilde D^3=5$. In the former
case the original distribution $D$ is locally equivalent to the
Goursat distribution and has an infinite-dimensional symmetry
algebra. In other words, the case of non-Goursat distributions of
constant class~1 can be reduced to the case of distributions of
class greater than~1.

This leaves the following question open: \emph{Do there exist
  completely nonholonomic rank $2$ distributions of constant class
  $2\leq m\leq n-4$?}
%between $2$ and $n-4$?
We know only that the answer is negative for $m=2$ ($n>5$),
which means that any such example, if it exists, should
live on at least $7$-dimensional manifold.

\subsection{Connection with Tanaka theory}

After the symplectification procedure described above, the
results of this paper can be interpreted in terms of
Tanaka--Morimoto theory of structures on filtered
manifolds~\cite{mori,tan2}. The original distribution $D$
(even of maximal class) has, in general, a non-constant
symbol, which makes this theory very difficult to apply to
the filtered manifold defined by the distribution $D$
itself. However, given rank 2 distribution $D$ of maximal
class there is a natural rank 2 distribution on the
manifold $P(\mathcal R_D)$ obtained from $\mathcal R_D$ via
the factorization by the trajectories of the Euler vector
field (or, in other words, by the projectivization of the
fibers of $\mathcal R_D$). It is generated by the
projection of the sum $V_{n-4}\oplus {\mathcal C}$ w.r.t.
this factorization. It is possible to show that this
distribution has already a fixed symbol isomorphic to the
Lie algebra generated by the vector fields
$\{h,\tilde\e_1,\dots,\tilde\e_{2m},\eta\}$ from the proof
of the main theorem (see equation~\eqref{commrel}).

Moreover, there is a natural decomposition of this
distribution into the sum of two line distributions equal
to the projections of $V_{n-4}$ and ${\mathcal C}$. This
decomposition can be interpreted as a $G$-structure on a
filtered manifold in terms of Tanaka theory and is called a
\emph{pseudo-product} structure~\cite{tan3}. The
prolongation of this structure (in terms of filtered
manifolds) is of finite type and is isomorphic to the
maximal symmetry algebra from the main theorem.

We shall dwell into the details of this approach in the forthcoming
paper.


\begin{thebibliography}{99}
%\bibitem{agrgam1}
%A.A. Agrachev, R.V. Gamkrelidze, Feedback-invariant optimal
%control theory - I. Regular extremals, J. Dynamical and
%Control Systems, {\bf 3}(1997), No. 3, 343-389.

\bibitem{agrachev}
A.A. Agrachev, Feedback-invariant optimal control theory - II.
Jacobi Curves for Singular Extremals, J. Dynamical and Control
Systems, {\bf 4}(1998), No. 4 , 583-604.

\bibitem{agrsach}
A.A. Agrachev, Yu.L. Sachkov Control Theory from the Geometric
Viewpoint Series: Encyclopaedia of Mathematical Sciences , Vol. 87,
412 pages.

\bibitem{agrzel}
A.A. Agrachev, I. Zelenko, Principal Invariants of Jacobi
Curves,
%In the book: Nonlinear Control in the Year 2000, v.
%1, A. Isidori, F. Lamnabhi-Lagarrigue $\&$ W. Respondek,
%Eds,
Lecture Notes in Control and Information Sciences 258,
Springer, 2001, 9-21.

\bibitem{jac1} A.A. Agrachev, I. Zelenko, Geometry of Jacobi
curves.I, J. Dynamical and Control Systems, {\bf 8},2002,
No. 1, 93-140.

\bibitem{jac2} A.A. Agrachev, I. Zelenko, Geometry of Jacobi
curves.II, J. Dynamical and Control Systems, {\bf 8},2002,
No.2, 167-215.


\bibitem{bryantbook}
R.L. Bryant, S.S. Chern, R. B. Gardner, H.L. Goldschmidt,
P. A. Griffiths, Exterior Differential Systems,
Mathematical Sciences Research Institute Publications, vol.
18, Springer-Verlag.

\bibitem{cartan}
E. Cartan, {\sl Les systemes de Pfaff a cinque variables et
les equations aux derivees partielles du second ordre},
%Ann. Sci. Ecole Normale 27(3), 1910, pp 109-192; reprinted
%in his
Oeuvres completes, Partie II, vol.2, Paris,
Gautier-Villars, 1953, 927-1010.

\bibitem{doubzel} B. Doubrov, I. Zelenko, {\sl A canonical frame for
nonholonomic rank two distributions of maximal class}, C.R. Acad.
Sci. Paris, Ser. I, Vol. 342, Issue 8 (15 April 2006), 589-594.

\bibitem{doubzel3}
B. Doubrov, I. Zelenko, {\sl Geometry of rank 3 distributions}, in
preparation.

\bibitem{mori} T. Morimoto, {\sl Geometric structures on
filtered manifolds}, Hokkaido Math. J.,{\bf 22}(1993), pp.
263-347.

\bibitem{tan2} N. Tanaka, {\sl On the equivalence problems associated with simple graded Lie
algebras}, Hokkaido Math. J.,{\bf 6}(1979), pp 23-84.

\bibitem{tan3} N.~Tanaka, {\sl On affine symmetric spaces and the
automorphism groups of product manifolds}, Hokkaido Math.
J., {\bf 14}(1985),
  pp 277-351.

\bibitem{versh1}
A. Vershik, V. Gershkovich, {\sl Determination of the
functional dimension of the orbit space of generic
distributions}, Mat. Zametki 44, 596-603 (in Russian);
English transl: Math. Notes 44, 806-810 (1988).

\bibitem{wil} E.J. Wilczynski, {\sl Projective Differential Geometry
of Curves and Ruled Surfaces}, Teubner, Leibzig, 1905.

\bibitem{zelvar} I. Zelenko, {\sl Variational Approach to Differential
Invariants of Rank 2 Vector Distributions}, Differential Geometry
and Its Applications, Vol. 24, Issue 3 (May 2006), 235-259
%to appear in
%"Differential Geometry and Its Applications", 30 pages;
% arxiv math. DG/0402171, SISSA preprint 12/2004/M.

\bibitem{praga} I. Zelenko, {\sl  Complete systems of invariants for rank 1 curves in Lagrange
Grassmannians}, Differential Geom. Application, Proc. Conf. Prague,
2005, pp 365-379, Charles University, Prague
%Proceedings of 9th
% Conference on Differential Geometry and its Applications,
% arxiv math. DG/0411190, SISSA preprint 82/2004/M, 15 pages

\bibitem {zel}
I. Zelenko, {\sl Nonregular abnormal extremals of
2-distribution: existence, second variation, and rigidity},
J. Dynamical and Control Systems, {\bf 5}(1999), No. 3,
347-383.

 \bibitem {zhit0} M. Zhitomirskii, {\sl Normal forms of germs of smooth
distributions}, Mat. Zametki 49 (1991), no. 2, 36--44, 158
(in Russian); English translation in Math. Notes 49 (1991),
no. 1-2, 139--144.

\end{thebibliography}
\end{document}